\documentclass[12pt, leqno]{article}

\usepackage{amssymb} 
\usepackage{amsmath}
\newtheorem{Thm}{Theorem}[section] 
\newtheorem{Lem}[Thm]{Lemma}
\newtheorem{Cor}[Thm]{Corollary} 
\newtheorem{Prop}[Thm]{Proposition}
\newenvironment{Pf}{\noindent {\it Proof.}\quad}{\quad $\square$\\ \\}
\newenvironment{Rmk}{\noindent {\it Remark.}}{\\}
\numberwithin{equation}{section}
\begin{document} 

\title {\bf Gr\"obner-Shirshov Bases for Lie
Superalgebras and Their Universal Enveloping Algebras
}

\author 
{Leonid A. Bokut \thanks{Supported in part by the Russian Fund of 
Basic Research.} \\ 
{\small Institute of Mathematics, Novosibirsk 630090, Russia} \\
{\small  School of Mathematics, Korea Institute for Advanced Study,
Seoul 130-010, Korea}
\and
Seok-Jin Kang \thanks{Supported in part by 
Research Institute of Mathematics
at Seoul National University and Korea Institute for Advanced Study.}
\\
{\small Department of Mathematics, Seoul National University,
Seoul 151-742, Korea}  \\ 
{\small School of Mathematics, Korea Institute for Advanced Study,
Seoul 130-010, Korea} 
\and
Kyu-Hwan Lee \thanks{Supported in part by Research Institute
of Mathematics and GARC-KOSEF at Seoul National University.}  \\
{\small
Department of Mathematics,
Seoul National University, Seoul 151-742, Korea}
\and
Peter Malcolmson  \\
{\small 
Department of Mathematics, Wayne State University,
Detroit, MI 48202, U.S.A. } }

\date{}

\maketitle

\begin{abstract}
We show that a set of monic polynomials in the free Lie superalgebra is a 
Gr\"obner-Shirshov basis for a Lie superalgebra if and only if
it is a Gr\"obner-Shirshov basis for its universal enveloping
algebra. We investigate the structure of 
Gr\"obner-Shirshov bases for Kac-Moody superalgebras and  
give explicit constructions of 
Gr\"obner-Shirshov bases for classical Lie superalgebras.

\end{abstract}

\section{Introduction} 

Let $\mathcal {A}$ be a free (commutative, associative, or Lie) 
algebra over a field $k$, let $S\subset \mathcal{A}$ be a set of relations
in $\mathcal {A}$, and let $\langle S \rangle$ be the ideal of 
$\mathcal {A}$ generated by $S$. 
One of the fundamental problems in the theory of abstract algebras is
the {\it reduction problem}: given an element $f \in \mathcal{A}$,
one would like to find a {\it reduced expression} for $f$ with respect to the 
relations in $S$.
One of the most common approaches to this problem 
is to find another set of generators for the relations in $S$ 
that can replace the original relations so that one can 
get an effective algorithm for the reduction problem. 
More precisely, if one can find a set $S^{c}$ of generators of the ideal
$\langle S \rangle$ which is {\it closed} under a certain {\it composition} of
relations in $S$, then there exists an easy criterion by which one can
determine whether an element $f \in \mathcal {A}$ is reduced 
with respect to $S$ or not.

In 1965, inspired by Gr\"obner's suggestion, 
Buchberger found a criterion and an algorithm of computing 
such a set of generators of the ideals for commutative algebras \cite{Bu1},
which were modified and refined in \cite{Bu2} and \cite{Bu3}.
Such a set of generators of ideals is now referred to as a 
{\it Gr\"obner basis}, and it has become one of the most popular 
research topics in the theory of commutative algebras (see for example,
\cite{BW}). 
In 1978, Bergman developed the theory 
of Gr\"obner bases for associative algebras by proving the
{\it Diamond Lemma} \cite{Be}. His idea is a generalization
of Buchberger's theory and it has many applications to various
areas of the theory of associative algebras such as quantum groups.

For the case of Lie algebras, where the situation is more complicated than
commutative or associative algebras, the parallel theory of Gr\"obner 
basis was developed by Shirshov in 1962 \cite{Sh3}, which is 
even earlier than Buchberger's discovery. In that paper, which 
was written in Russian and never translated in English,  
he introduced the notion of 
{\it composition} of elements of a free Lie algebra
and showed that a set of relations which is closed under the composition 
has the desired property.  
Shirshov's idea is essentially the same as that of Buchberger, 
and it was noticed by Bokut that Shirshov's method works for associative 
algebras as well \cite{Bo3}. 
For this reason, we will call such a set of relations of a free Lie algebra 
(and of a free associative algebra) a {\it Gr\"obner-Shirshov basis}.
(See \cite{BMM} for a more detailed history
of Gr\"obner-Shirshov basis.)
It has been used to determine the solvability of some 
word problems \cite{Sh2, Sh3, Bo2} and to prove some embedding theorems 
\cite{Bo1, Bo3, Bo4}. 
Recently, in a series of works by Bokut, Klein, and Malcolmson, 
Gr\"obner-Shirshov bases for finite dimensional 
simple Lie algebras and the quantized enveloping algebra of type $A_n$ 
were constructed explicitly (\cite{BoKl1, BoKl2, BoKl3, BM}).

In this work, we develop the theory of Gr\"obner-Shirshov bases for
Lie superalgebras and their universal enveloping algebras. 
This paper is organized as follows. 
In Section 2, after introducing the basic facts
such as {\it super-Lyndon-Shirshov words} ({\it monomials}) 
and {\it Composition
Lemma}, we prove that a set of monic polynomials in the free
Lie superalgebra is a 
Gr\"obner-Shirshov basis for a Lie superalgebra if and only if
it is a Gr\"obner-Shirshov basis for its universal enveloping
algebra (Theorem \ref{thm:2main}). 
This is a generalization of the corresponding result for 
Lie algebras obtained in \cite{BoM}. 
Thus the theory of Gr\"obner-Shirshov bases for Lie superalgebras
and that of associative algebras are unified in this way, and
as a by-product, we obtain a purely combinatorial proof of 
the Poincar\'e-Birkhoff-Witt Theorem (Proposition \ref{prop:PBW}).

In section 3, we investigate the structure of Gr\"obner-Shirshov
bases for Kac-Moody superalgebras and prove that, in order to 
find a Gr\"obner-Shirshov basis for a Kac-Moody superalgebra,
it suffices to consider the completion of Serre relations of the
positive part (or negative part) which is closed under the composition 
(Theorem \ref{thm-a}). 
As a corollary, we obtain the {\it triangular decomposition} of
Kac-Moody superalgebras and their universal enveloping algebras 
(Corollary \ref{cor-a}). 
Our result in this section is a generalization of the corresponding 
result for Kac-Moody algebras obtained in \cite{BM}.

Finally, in Section 4, we give an explicit construction of
Gr\"obner-Shirshov bases for classical Lie superalgebras. 
The outline of our construction can be described as follows. 
We first start with a Kac-Moody superalgebra which is
isomorphic to a given classical Lie superalgebra. 
Using the supersymmetry and Jacobi identity, we expand the
set of Serre relations to a complete set $R$ of relations 
which is closed under the composition and
determine the set $B$ of $R$-reduced super-Lyndon-Shirshov monomials.
Now comparing the number of elements of $B$ with the 
dimension of the corresponding classical Lie superalgebra, we conclude that
the set $R$ is indeed a Gr\"obner-Shirshov basis. 

The main part of this work was completed while the first and the second
authors were visiting Korea Institute for Advanced Study in the winter
of 1998. We would like to express our sincere gratitude to
Hyo-Chul Myung, Efim Zelmanov and the other members of 
Korea Institute for Advanced Study for their hospitality
and cooperation.

\vskip 1cm

\section{Gr\"obner-Shirshov bases for Lie superalgebras} 

Let $X=X_{\bar{0}} \cup X_{\bar{1}}$ be a
$\mathbb{Z}_{2}$-graded set with a linear ordering $\prec$, 
and let $X^*$ (resp. $X^{\#}$) be 
the semigroup of associative words on $X$ (resp. the groupoid of
nonassociative words on $X$).  
Then the semigroup $X^*$ (resp. the groupoid $X^\#$)
has the $\mathbb{Z}_{2}$-grading $X^*=X^*_{\bar{0}} \oplus
X^*_{\bar{1}}$ (resp.  $X^\#=X^\#_{\bar{0}} \oplus X^\#_{\bar{1}}$) induced by
that of $X$.  
The elements of $X^*_{\bar{0}}$ and $X^{\#}_{\bar 0}$ (resp. 
$X^*_{\bar{1}}$ and $X^{\#}_{\bar 1}$) 
are called {\it even} (resp. {\it odd}).  

We denote by $l(u)$ the {\it length} of a word $u$
and the empty word will be denoted by  $1$. 
For an associative word $u \in X^*$, we can choose a certain 
arrangement of brackets on $u$, which will be denoted by $(u)$. 
Conversely, there is a canonical bracket removing homomorphism 
$\rho :  X^\# \rightarrow X^*$ given by 
$\rho( (u) )=u$ for $u \in X^*$.  

We consider two linear orderings $<$ and $\ll$ 
on $X^*$ defined as follows:  
\begin{list}{ }
{\setlength{\itemsep}{0cm} \setlength{\leftmargin}{1cm}} 
\item (i) $u<1$ for any
nonempty word $u$; and inductively, $u<v$ whenever $u=x_i u'$, $v=x_j v'$ and
$x_i \prec x_j$ or $ x_i=x_j$ and $u' < v'$.  
\item (ii) $u \ll v$ if $l(u) <
l(v)$ or $l(u) = l(v)$ and $u<v$.  \end{list} 
The ordering $<$ (resp. $\ll$) is called the {\it lexicographical 
ordering} (resp. {\it length-lexicographical ordering}).
We define the orderings $<$ and $\ll$ on $X^\#$ by 
(i) $u<v$ if and only if 
$\rho(u) < \rho(v)$, and (ii) $u \ll v$ if
and only if $\rho(u) \ll \rho(v)$.  

A nonempty word $u$ is called a {\it Lyndon-Shirshov word} 
if $u \in X$ or $vw > wv$ for any
decomposition of $u=vw$ with $v, w \in X^*$.
A nonempty word $u$ is called a 
{\it super-Lyndon-Shirshov word} if either it is a Lyndon-Shirshov word 
or it has the form $u=vv$ with $v$ a Lyndon-Shirshov 
word in $X^*_{\bar{1}}$.  
A nonempty nonassociative word $u$ is called a 
{\it Lyndon-Shirshov monomial} if either $u$ is an element of $X$ 
or  
\begin{list}{ } {\setlength{\itemsep}{0cm} \setlength{\leftmargin}{1cm}}
\item (i) if $u=u_1 u_2$, then $u_1, u_2$ are 
 Lyndon-Shirshov monomials with $u_1 > u_2$, 
\item (ii) if $u=(v_1 v_2)w$ then $v_2 \leq w$.  \end{list} 
A nonempty nonassociative word $u$ is called a
{\it super-Lyndon-Shirshov monomial} if either it is 
a Lyndon-Shirshov monomial 
or it has the form $u=vv$ with $v$ 
a Lyndon-Shirshov monomial in $X^\#_{\bar{1}}$.

\begin{Rmk}
In some literatures, the Lyndon-Shirshov words have been referred to
as {\it regular words}, {\it normal words}, {\it Lyndon words}, etc. 
Since the definition of Lyndon-Shirshov words dates back to the works by
Chen, Fox and Lyndon \cite{CFL} and Shirshov \cite{Sh1}, we decide
to call them Lyndon-Shirshov words. 
The definition of super-Lyndon-Shirshov words can be found in
\cite{Ba, Mik}.
\end{Rmk}

The following lemma asserts that there is a natural 
1-1 correspondence between the set of super-Lyndon-Shirshov
words and the set of super-Lyndon-Shirshov monomials. 

\begin{Lem} {\rm (\cite{Ba, CFL, Mik, Sh0})}
If $u$ is a super-Lyndon-Shirshov monomial, 
then $\rho(u)$ is a super-Lyndon-Shirshov word.
Conversely, for any super-Lyndon-Shirshov word $u$,
there is a unique arrangement of brackets
$[u]$ on $u$ such that $[u]$ is a super-Lyndon-Shirshov monomial.  
\end{Lem}

Let $k$ be a field with $\mathrm{char}(k) \neq 2,3$, 
and let $\mathcal{A}_{X}$ be the
free associative algebra generated by $X$ over $k$.  The algebra
$\mathcal{A}_{X}$ becomes a Lie superalgebra with the
superbracket defined by 
$$[x,y] = xy-(-1) ^{(\deg x) (\deg y)} yx$$ 
for $x, y \in \mathcal{A}_{X}$.  
Let $\mathcal{L}_{X}$ be the subalgebra of $\mathcal{A}_{X}$ 
generated by $X$ as a Lie superalgebra.  
Then $\mathcal{L}_{X}$ is the free Lie superalgebra generated by
$X$ over $k$.  
As we can see in the following theorem, there is a 
canonical linear basis for the free Lie superalgebra $\mathcal{L}_{X}$:

\begin{Thm} {\rm (\cite{Ba, CFL, Mik, Sh0})}
The set of super-Lyndon-Shirshov monomials form a linear basis of 
the free Lie superalgebra $\mathcal{L}_{X}$ generated by $X$.
\end{Thm} 

\begin{Rmk} The existence of linear bases for free Lie algebras of this form 
was first suggested by Hall \cite{Hall}, and later by Shirshov
in a more general form (\cite{Sh0, Sh4}).  
The linear basis for a free  Lie superalgebra
given in the above theorem will be called the
{\it Lyndon-Shirshov basis}.  It is a special case of the 
{\it Hall-Shirshov basis}. \end{Rmk}

Given a nonzero element $p \in \mathcal{A}_{X}$ 
we denote by $\overline{p}$ the
maximal monomial appearing in $p$ under the ordering $\ll$.  
Thus $p = \alpha
\overline{p} + \sum \beta _i w_i $ with $\alpha , \beta _i \in k$, $ w_i \in
X^*$, $\alpha \neq 0$ and $w_i \ll \overline{p}$.  
The coefficient $\alpha$ of $\overline{p}$ is called the
{\it leading coefficient} of $p$ and $p$ is said to be {\it
monic} if $\alpha =1$.  

The following lemma plays a crucial role in defining 
the notion of {\it Lie composition}.

\begin{Lem} \label{lem-2ab} {\rm (\cite{CFL, Mik, Sh0})} 
Let $u$ and $v$ be super-Lyndon-Shirshov words such that 
$v$ is contained in $u$ as a subword. Write 
$u=avb$ with $a, b \in X^*$.  
Then there is an 
arrangement of brackets $[u] = (a[v]b)$ on $u$ such that $[v]$ is a
super-Lyndon-Shirshov monomial, $\overline{[u]}=u$ and 
the leading coefficient of $[u]$ is either $1$ or $2$.
\end{Lem}

Let $u=avb$ be a super-Lyndon-Shirshov word, where 
$v$ is a super-Lyndon-Shirshov subword and $a, b \in X^*$.
We define the {\it bracket on $u$ relative to $v$}, 
denoted by $[u]_v$, as follows:  
\begin{list}{ } 
{\setlength{\itemsep}{0cm} \setlength{\leftmargin}{1cm}}
\item (i) \, $[u]_v =  (a[v]b)$ if the leading coefficient of $[u]$ is $1$, 
\item (ii) \, $[u]_v =  \frac{1}{2} (a[v]b)$ if the leading coefficient
of $[u]$ is $2$, 
\end{list} 
where the arrangement of brackets $[u]$ on $u$ is the one 
described in Lemma \ref{lem-2ab}.  
Note that $[u]_v$ is monic and $\overline{[u]_v} =u$.

Similarly, if $p$ is a monic polynomial in the free Lie superalgebra 
$\mathcal{L}_{X}$ such that $\overline {p}$ is super-Lyndon-Shirshov, 
then we define the {\it bracket on $u$ relative to $p$},
denoted by $[u]_{p}$ to be the result of the substitution of $p$ 
instead of $\overline {p}$ in $[u]_{\overline {p}}$. 
Clearly, $[u]_{p}$ is monic and $\overline{[u]_{p}}=u$. 

We now define the notion of {\it associative composition} 
of the elements in the free associative algebra $\mathcal{A}_{X}$
generated by $X$. 
Let $p, q$ be monic elements in $\mathcal{A}_{X}$ 
with leading terms $\overline{p}$ and $\overline{q}$.  
If there exist $a, b \in X^*$ such that $\overline{p}a =
b \overline{q} = w$ with $l(\overline{p}) > l(b)$, then
we define the {\it composition of
intersection} $(p,q)_{w}$ to be
\begin{equation}
(p,q)_w = pa -bq.
\end{equation}
If there exist $a, b \in X^*$ such that
$\overline{p}=a\overline{q}b=w$, then we define
the {\it composition of inclusion} to be
\begin{equation} 
(p,q)_w = p - aqb.
\end{equation}
Note that we have $\overline{(p,q)_{w}} \ll w$ 
in either case.

Next we proceed to define the notion of 
{\it Lie composition} of the elements in the free Lie superalgebra
$\mathcal{L}_{X}$ generated by $X$. 
Let $p$, $q$ be monic polynomials in the free Lie superalgebra 
$\mathcal{L}_{X}$ with leading terms $\overline{p}$
and $\overline{q}$.  If there exist $a, b \in X^*$ such that
$\overline{p}a = b \overline{q} = w$ with $l(\overline{p}) > l(b)$, 
then we define the {\it composition of intersection} 
$\langle f,g \rangle _w$ to be
\begin{equation}
\langle f,g \rangle _w = [w]_p -[w]_q.
\end{equation}  
If there exist 
$a, b \in X^*$ such that $\overline{p}=a\overline{q}b=w$, then
we define the {\it composition of inclusion} to be
\begin{equation}
\langle p,q \rangle _w = p - [w]_q.
\end{equation}
We have 
$\overline{\langle p, q \rangle _w} \ll w$
in this case, too.   \\

\begin{Rmk} Our definition of Lie composition is essentially the same as
the one given in \cite{Bo2, Lo, Mik, Sh2}.  
We modified the definition in \cite{Bo2, Lo, Mik, Sh2} 
to define the Lie composition $\langle p, q \rangle_w$ 
at one stroke.  \end{Rmk}

Let $S$ be a set of monic polynomials in $\mathcal{L}_{X} \subset
\mathcal{A}_{X}$, let $I$ be
the (Lie) ideal generated by $S$ in the free Lie superalgebra 
$\mathcal{L}_{X}$, and let $J$ be the (associative) ideal generated
by $S$ in the free associative algebra $\mathcal{A}_{X}$.  
We denote by $L = \mathcal{L}_{X} / I$  the Lie
superalgebra generated by $X$ with defining relations $S$
and let $\mathcal{U}(L) = \mathcal{A}_{X}/J $ be its universal 
enveloping algebra.  

For $f, g \in \mathcal{A}_{X}$ and $w \in X^*$, we write $f \equiv_A g$
$ \mathrm{mod} \, (S, w)$ if $f -g = \sum \alpha_i a_i s_i b_i$, where
$\alpha_i \in k, a_i,b_i \in X^*, s_i \in S$ 
with $a_i \overline{s_i} b_i \ll w$
for each $i$. 
Similarly, for $f, g \in \mathcal{L}_{X}$ and $w \in X^*$, 
we write $f \equiv_L g$
$ \mathrm{mod} \, (S, w)$ if $f -g = \sum \alpha_i (a_i (s_i) b_i)$, where
$\alpha_i \in k, a_i,b_i \in X^*, s_i \in S$ with 
$\overline{p(a_i (s_i) b_i)} \ll w$
for each $i$. The set $S$ is said to be {\it closed under the associative 
composition} (resp. {\it Lie composition}) 
if for any $f, g \in S$, we have $(f,g)_w \equiv_A 0$
(resp.  $\langle f,g \rangle _w \equiv_L 0$) $\mathrm{mod} \, (S, w)$.  

A set of monic polynomials $S$ in the free Lie superalgebra 
$\mathcal{L}_{X}$ is called a {\it Gr\"obner-Shirshov basis} for the ideal $J$
(resp. for the ideal $I$) 
if it is closed under the associative composition
(resp. Lie composition).
By abuse of language, we will also
refer to $S$ as a Gr\"obner-Shirshov basis for the associative algebra 
$\mathcal{U}(L)$ and for the Lie superalgebra $L$, respectively. 
An associative word $u$ is said to be {\it S-reduced} if $u \neq
a\overline{s}b$ for any $s \in S$ and $a, b \in X^*$.  
A nonassociative word $u$
is said to be {\it S-reduced} if $\rho(u)$ is S-reduced.  

The following lemma is a generalization of Lemma 1 in \cite{BoKl1}.

\begin{Lem} \label{lem-aa} 
\hfill

{\rm (a)} Every nonempty word $u$ in the free associative algebra 
$\mathcal{A}_{X}$ can be written as 
\begin{equation}
 u = \sum \alpha_i u_i + \sum
\beta_j a_js_jb_j , 
\end{equation} 
where $u_i$ is an S-reduced word, $\alpha_i, \beta_j \in k$, 
$a_j, b_j \in X^*$, $s_j \in S$ and $a_j\overline{s_j}b_j \,
\underline{\ll} \, u$ for all $i, j$.
Hence the set of $S$-reduced words spans the algebra $\mathcal{U}(L)$.

{\rm (b)} Every super-Lyndon-Shirshov monomial $u$ in
$\mathcal{L}_{X}$ can be written as
\begin{equation}
u = \sum \alpha_i u_i + \sum \beta_j
(a_j (s_j) b_j), 
\end{equation} 
where $u_i$ is an S-reduced super-Lyndon-Shirshov monomial, 
$\alpha_i, \beta_j \in k$, $a_j, b_j \in X^*$, $s_j \in S$ 
and $\overline{(a_j (s_j) b_j)} \,
\underline{\ll} \, \overline{u}$ for all $i,j$.
Hence the set of S-reduced 
super-Lyndon-Shirshov monomials spans the Lie superalgebra $L$. 
\end{Lem}

\begin{Pf} Since the proof of (a) is similar to that of (b),
we only give a proof of (b). 
If $u$ is S-reduced, we are done.  
Thus we assume that $\overline{u}
=a\overline{s}b$ for some $s \in S$, $a, b \in X^*$.  
Then $\overline{u}$ and $\overline{s}$ are
super-Lyndon-Shirshov words and 
$\overline{u-\alpha[\overline{u}]_s} \ll \overline{u}$ for some
$\alpha \in k$.  Since $u-\alpha[\overline{u}]_s$ is a linear combination of
super-Lyndon-Shirshov monomials whose 
leading terms are less than $\overline{u}$, we may proceed by
induction, which completes the proof.
\end{Pf}

The following lemma plays a crucial role in our discussion of
Gr\"obner-Shirshov bases. It is originally due to Shirshov \cite{Sh3}
and is now known as the {\it Composition Lemma}.

\begin{Lem} \label{lem-comp} {\rm (cf. \cite{Ba, Bo2, Mik, Sh3})}
If $S$ is a Gr\"obner-Shirshov basis for the ideal $J$, 
then for any $f \in J$, the word
$\overline{f}$ contains a subword $\overline{s}$ with $s \in S$.  
\end{Lem}

It is clear that if a polynomial $f\in \mathcal{L}_{X}$ satisfies 
$f \equiv_L 0$ $\mathrm{mod}(S, w)$ for $w\in X^*$, then 
$f \equiv_A 0$ $\mathrm{mod}(S, w)$. 
The converse is also true if $S$ is closed under the associative 
composition. 

\begin{Lem} \label{lem-cong}
Assume that $S$ is closed under the associative composition. 
If a polynomial $f \in \mathcal{L}_X$ satisfies 
$f \equiv_A 0$ $\mathrm{mod}(S, w)$ for $w \in X^*$, 
then $f \equiv_L 0$ $\mathrm{mod}(S, w)$.
\end{Lem}

\begin{Pf}
Suppose $f \equiv_A 0$ $\mathrm{mod}(S, w)$ for $w\in X^*$ and 
our assertion holds for all $w' \ll w$. 
Then $f \in J$, and by the Composition Lemma, 
$\overline{f} = a\overline{s} b$ for some $a, b\in X^*$ and $s \in S$. 
Since $f - [\overline{f}]_{s} \equiv_A 0$ $\mathrm{mod}(S, \overline{f})$ 
and $\overline{f} \ll w$, our assertion follows by induction.
\end{Pf}

\begin{Lem}
\label{lem-equi} Let $f, g \in S$ be monic polynomials in
$\mathcal{L}_{X}$ such that the associative composition $(f, g)_w$
is defined.  
Then we have 
\begin{equation}
 (f, g)_w \equiv_A \langle f, g \rangle_w \ \ \ \mod (S, w).
\end{equation}
\end{Lem}

\begin{Pf} We consider the composition of intersection only.
The proof for the composition of inclusion is similar. 
Recall that $[w]_f = fa + \sum
\alpha_i a_ifb_i$ with $a_i\overline{f}b_i \ll w$ 
and $[w]_g = bg + \sum \beta_i c_igd_i$ 
with $c_i\overline{g}d_i \ll w$.  
Thus $\langle f, g \rangle_w = [w]_f
- [w]_g = fa-bg +h = (f,g)_w +h$, 
where $h \equiv_A 0$ $\mathrm{mod} (S, w)$.  
Hence $(f,g)_w \equiv_A \langle f, g \rangle _w$ $\mathrm{mod} (S, w)$. 
\end{Pf}

Combining Lemma \ref{lem-cong} and Lemma \ref{lem-equi},
we obtain the main result of this section, which is a 
generalization of the main theorem in \cite{BoM}. 

\begin{Thm} \label{thm:2main}
Let $S$ be a set of monic polynomials in the free Lie superalgebra
$\mathcal{L}_{X}$. Then $S$ is a Gr\"obner-Shirshov basis 
for the Lie superalgebra $L=\mathcal{L}_{X} \big/ I$ 
if and only if $S$ is a Gr\"obner-Shirshov basis for 
its universal enveloping algebra $\mathcal{U}(L)=\mathcal{A}_{X} \big/J$.
That is, $S$ is closed under the Lie composition if and only if it 
is closed under the associative composition. 
\end{Thm}

The following proposition, which is a generalization of 
Proposition 2 in \cite{BoKl1}, provides us with a criterion 
for determining whether a set of monic polynomials in the free Lie
superalgebra is a Gr\"obner-Shirshov basis or not.  

\begin{Prop} \label{prop-tt} 
\hfill

{\rm (a)} If the set of S-reduced words is a linear basis of
$\mathcal{U}(L)=\mathcal{A}_{X}\big/ J$, 
then $S$ is a Gr\"obner-Shirshov basis for the
ideal $J$ of $\mathcal{A}_{X}$.

{\rm (b)} If the set of S-reduced super-Lyndon-Shirshov monomials 
is a linear basis of $L=\mathcal{L}_{X} \big/ I$, then $S$ is a
Gr\"obner-Shirshov basis for the ideal $I$ of $\mathcal{L}_{X}$.
\end{Prop}

\begin{Pf} Since the proof of (b) is the same as (a), we will prove
(a) only.
Suppose on the contrary that $S$ is not closed under
the associative composition.
Then there exist $f, g \in S$ such that $(f,g)_w \not\equiv_{A} 0$ 
$\mathrm{mod} (S, w)$ for
$w\in X^*$. By Lemma \ref{lem-aa}, we may write 
$$ (f,g)_w =
\sum \alpha_i u_i + \sum \beta_j a_js_jb_j , $$ 
where $\alpha_i, \beta_j \in k$,
$u_i$ is S-reduced, $a_j, b_j \in X^*$, $s_j \in S$ and $a_j\overline{s_j}b_j
\ll w$ for all $i$ and $j$.  
Since $(f,g)_w \not\equiv_{A} 0$ $\mathrm{mod} (S, w)$, 
we have $\sum \alpha_i u_i \neq 0$ in $\mathcal{A}_{X}$.
Since the set of S-reduced words is a linear basis of $\mathcal{U}(L)$,
we have $\sum \alpha_i u_i \neq 0$ in $\mathcal{U}(L)$.  
But, since $(f,g)_w \in J$, we have $\sum \alpha_i u_i
= 0$ in $\mathcal{U}(L)$, which is a contradiction. 
\end{Pf}

Conversely, by Lemma \ref{lem-aa} and the Composition Lemma, 
we can show that a Gr\"obner-Shirshov basis gives rise to a
linear basis for the corresponding algebras. 

\begin{Thm} \label{thm-b} 
\hfill 

{\rm (a)} If $S$ is a Gr\"obner-Shirshov basis for the 
Lie superalgebra $L = \mathcal{L}_{X} / I$, then
the set of S-reduced super-Lyndon-Shirshov monomials forms a 
linear basis of $L$. 

{\rm (b)} If $S$ is a Gr\"obner-Shirshov basis for the universal enveloping
algebra $\mathcal{U}(L) = \mathcal{A}_{X}/J$ of $L$,
then the set of S-reduced words forms a linear basis of 
$\mathcal{U}(L)$.  
\end{Thm}

\begin{Pf}
Since the proof of (b) is similar to that of (a), we will prove (a) only. 
 By Lemma \ref{lem-aa} the set of $S$-reduced super-Lyndon-Shirshov monomials
spans $L$. Assume that we have $ \sum \alpha_i u_i =0$ in $L$, where 
$\alpha_i \in k$ and  $u_i$ are distinct $S$-reduced super-Lyndon-Shirshov 
monomials. Then $\sum \alpha_i u_i \in I$ in the free Lie super algebra 
$\mathcal{L}_X$. Since $I \subset J$, we obtain
$\sum \alpha_i u_i \in J$. By the Composition Lemma (Lemma \ref{lem-comp})
the leading term $\overline{ \sum \alpha_i u_i}$ contains a subword 
$\overline{s}$ with 
$s \in S$. Since each $u_i$ is $S$-reduced, we must have $\alpha_i = 0$
for all $i$. Hence the set of $S$-reduced super-Lyndon-Shirshov monomials
is linearly independent.
\end{Pf}

As a corollary, we obtain a purely combinatorial proof of the
Poincar\'e-Birkhoff-Witt Theorem. 

\begin{Prop} \label{prop:PBW}
\hfill

Let $L = L_{\bar{0}} \oplus L_{\bar{1}}$ be a Lie superalgebra with a linear
basis $Z= \{ z_1, z_2, \dots \}$ such that each $z_i$ is homogeneous with 
respect to the $\mathbb{Z}_{2}$-grading.
Then a linear basis of the universal enveloping algebra $\mathcal{U}(L)$
of $L$ is given by the set of all elements of the form $ z_{i_1} z_{i_2}
\dots z_{i_n}$ where $i_{k} \le i_{k+1}$ and $i_{k} \neq i_{k+1}$ if $z_{i_k}
\in L_{\bar{1}}$.
\end{Prop}

\begin{Pf}
Let $Y=\{ y_1, y_2, \dots \}$ be a $\mathbb{Z}_{2}$-graded set 
identified with the set $Z$ by a map $\iota$ such that $\iota (y_i) = z_i$ 
and $\iota (Y_{\alpha}) = Z_{\alpha}$ with $\alpha \in \mathbb{Z}_{2}$. 
Let $\mathcal{L}_Y$ be the free Lie superalgebra generated by 
$Y$. Let $S \subset \mathcal{L}_Y$ be the set of elements of the form 
$$ [y_i y_j] - \sum_k \alpha^k _{ij} y_k$$
where $i \ge j$ and  $i \neq j$ if $y_i \in Y_{\bar{0}}$,
and  $\alpha^k_{ij}$ is the structure constants given by the equation
$ [z_i z_j] = \sum_k \alpha^k _{ij} z_k $ in $L$. 
Let $I$ be the ideal of $\mathcal{L}_Y$ generated by $S$. Then, clearly, 
$\mathcal{L}_Y \big / I$ is isomorphic to $L$ and the set of $S$-reduced 
super-Lyndon-Shirshov monomials is just the set $Y$. By Proposition
\ref{prop-tt} the set $S$ is a Gr\"obner-Shirshov basis for $L$ and then 
by Theorem \ref{thm:2main} the set $S$ is also a Gr\"obner-Shirshov basis 
for $\mathcal{U}(L)$. Now our assertion follows from Theorem \ref{thm-b}. 
\end{Pf}

Let $S$ be a set of relations in the free Lie superalgebra $\mathcal{L}_{X}$
generated by $X$. We will see how one can {\it complete} the set $S$
to get a Gr\"obner-Shirshov basis. For any subset $T$ of $\mathcal{L}_{X}$,
we define 
$\widehat{T} = \{ p/\alpha \, | \,
\mbox{$\alpha \in k$ is the leading coefficient of $p \in T$} \}.$
Let $S^{(0)} = \widehat {S}$ and 
$S_{(0)} = \{ \langle f, g \rangle_w \not\equiv_{L} 0 \ 
\mbox{mod}\, (S^{(0)}, w) | \, f, g \in S^{(0)} \}. $
For $i\ge 1$, set 
$S_{(i)} = \{ \langle f, g \rangle_w \not\equiv_{L} 0 \ 
\mbox{mod}\, (S^{(i)}, w) | \, f, g \in S^{(i)} \}$
and 
$S^{(i)} =S^{(i-1)} \cup \widehat {S}_{(i-1)}.$

Then the set $S^c = \bigcup_{i\ge 0} S^{(i)}$
is a Gr\"obner-Shirshov basis for the (Lie) ideal $I$ 
generated by $S$ in $\mathcal{L}_X$.
Hence,  by Lemma \ref{lem-equi}, it is also a Gr\"obner-Shirshov 
basis for the (associative) ideal $J$ generated by $S$ in $\mathcal{A}_X$.  
It is easy to see that if every element of $S$ is homogeneous 
in $x_i \in X$,
then every element of $S^c$ is also homogeneous in $x_i$'s.

\section{Kac-Moody superalgebras} 

We now investigate the structure of Gr\"obner-Shirshov bases for 
Kac-Moody superalgebras. Our result is a generalization of the
work by Bokut and Malcolmson \cite{BM} on the Gr\"obner-Shirshov bases
for Kac-Moody algebras.
In the section, since we will consider the associative congruences only,
we will use the notation $\equiv$ in place of $\equiv_A$.

Let $\Omega=\{1, 2, \dots, r\}$ be a finite index set and $\tau$ be a subset
of $\Omega$. A square matrix $A=(a _{ij})_{i, j \in \Omega}$ is called
a {\it generalized Cartan Matrix} if it satisfies:
$$ \begin{array}{l} 
\mbox{(i) $a_{ii} =2$ or $0$ for $i= 1, \dots, r$ and if
$a_{ii} = 0$, then $i \in \tau$,}\\ 
\mbox{(ii) if $a_{ii} \neq 0$, then $a_{ij}
\in \mathbb{Z}_{\le 0}$ for $i \neq j$,}\\ 
\mbox{(iii) $a_{ij} =0$ implies $a_{ji} =0$,}\\ 
\mbox{(iv) if $a_{ii}=2$ and $i \in \tau$, then $a_{ij} \in 2
\mathbb{Z}$.}  \end{array} $$ 
Let $E=\{e_i\}_{i \in \Omega} ,
H=\{h_i\}_{i \in \Omega} , F=\{f_i\}_{i \in \Omega}$, 
and $X= E \,\cup H \, \cup F$.  
We define a $\mathbb{Z}_{2}$-grading on $\Omega$ by setting  $\deg i =
\bar{0}$ for $i \notin \tau$ and $\deg i = \bar{1}$ for $i \in \tau$, 
and on $X$ by  $\deg e_i =  \deg f_i = \deg i$ and $\deg h_i = \bar{0}$.  
We give a linear ordering on $X$ by  $e_{i} \succ h_{j} \succ f_{k}$ 
for all $i, j, k \in \Omega$ and $e_i \succ e_j$, $h_i \succ h_j$,
$f_i \succ f_j$ when $i > j$. Then we have the lexicographic ordering
and length-lexicographic ordering as in Section 2.  
We denote the left adjoint action of a Lie algebra by $\mathrm{ad}$
and the right adjoint action by $\widetilde{\mathrm{ad}}$.  The {\it Kac-Moody
superalgebra} $\mathcal{G} = \mathcal{G}(A, \tau)$ {\it associated to} $(A,
\tau)$ is defined to be the Lie superalgebra with generators $X$ and the
following defining relations:\\ 
\begin{equation}
\begin{aligned}
W:  & \quad [h_i h_j] \quad
(i>j),\\ & \quad [e_if_j] - \delta_{ij}h_i, \quad [e_jh_i] + a_{ij}e_j,\quad
[h_i f_j] + a_{ij}f_j, \\ 
S_{+,1}:  & \quad (\mathrm{ad}e_i)^{1-n_{ij}}e_j \quad
(i>j),\\ & \quad e_i (\widetilde{\mathrm{ad}}e_j)^{1-n_{ji}} \quad (i>j), \\ 
S_{+,2}:
& \quad [[e_{k+1}, e_k][e_k,e_{k-1}]] \quad \mbox{for $k \in \eta$},\\ 
S_{-,1}:
& \quad (\mathrm{ad}f_i)^{1-n_{ij}}f_j \quad (i>j),\\ & \quad f_i
(\widetilde{\mathrm{ad}}f_j)^{1-n_{ji}} \quad (i>j), \\ 
S_{-,2}:  & \quad [[f_{k+1},
f_k][f_k,f_{k-1}]] \quad \mbox{for $k \in \eta$}, 
\end{aligned}
\end{equation}
where 
\begin{equation}
n_{ij} = \left\{ \begin{array}{cl} 
a_{ij} & \mbox{if $a_{ii} = 2$ or $a_{ij} = 0$} \\ 
-1 & \mbox{if $a_{ii} = 0$ and $a_{ij} \neq 0$} \end{array} \right.  
\quad \mbox{for $i \neq j$}
\end{equation}
and $\eta$ is the set of indices $k$ such that $k \in \tau$,
$k \pm 1 \notin \tau$, $a_{kk}=0$, $a_{k+1,k-1}=0$ and $a_{k,k+1} +
a_{k,k-1}=0$.  Let $S_{\pm} = S_{\pm,1} \cup S_{\pm, 2}$ and $S(A, \tau) = S_+
\cup W \cup S_-$.  We denote by $\mathcal{G}_+$ (resp.  $\mathcal{G}_0$ and
$\mathcal{G}_-$) the subalgebra of $\mathcal{G}$ generated by $E$ (resp.  $H$
and $F$).

Set $t_{ij} =[e_if_j] - \delta_{i,j}h_i$, which belong to the relations $W$.  
We define the {\it differential substitution}
$\tilde{\partial}_j=\tilde{\partial}(e_j \rightarrow h_j)$ acting as a right
superderivation on $\mathcal{A}_E$ by  
\begin{equation}
\begin{aligned}
(e_i)\tilde{\partial}_j &= \delta_{ij}h_j, \\ 
(uv)\tilde{\partial}_j & =
u(v)\tilde{\partial}_j+(-1)^{(\deg j)(\deg v)}(u)\tilde{\partial}_j \, v \quad
\mbox{for $u, v \in \mathcal{A}_E$.} 
\end{aligned}
\end{equation}
It is easy to prove that for any $p \in \mathcal{A}_E$,
\begin{equation} \label{eq-a} pf_j \equiv (-1)^{(\deg p)(\deg j)} f_j p + (p)
\tilde{\partial}_j \quad \mod ( W,w).  
\end{equation} for some $w \gg
\overline{p}f_j$.  Note that $\tilde{\partial}_j$ is also a right
superderivation on $\mathcal{L}_E$.

\begin{Lem} \label{lem-a} 
Let $p$ be a homogeneous monic element of
$\mathcal{A}_E$ such that $(p, t_{ij})_w$ is defined for $w \in X^*$.  
Then we have 
$$(p, t_{ij})_w \equiv (p)\tilde{\partial}_j 
\quad \mod (\{p\} \cup W, w).  $$ 
\end{Lem}

\begin{Pf} 
It suffices to consider the composition of intersection.  We can write
$p=\overline{p} +p'$ with $\overline{p}=be_i$, where 
all the terms of $p'$ are lower than $\overline{p}$.
Then $w=\overline{p}f_j=be_if_j$.  Since $p$ is
homogeneous, $\deg p = \deg p'$.  From (\ref{eq-a}), we have 
\begin{equation*}
\begin{aligned}
 (p, t_{ i j } )_w & = pf_j - b ( e_i f_j - (-1)^{ (\deg i)
(\deg j ) }f_j e_i - \delta_{ij} h_j)\\ 
& = p'f_j + (-1)^{(\deg i)(\deg j)}
bf_je_i+\delta_{ij}bh_j\\ 
& \equiv (-1)^{ (\deg p) (\deg j ) } (f_j p' + f_jbe_i
) + (p') \tilde{\partial}_j \\
& +(-1)^{ (\deg i)(\deg j ) } ( b)\tilde{\partial}_j
e_i + \delta_{ij}b h_j\\ 
& \equiv (-1)^{ (\deg p) (\deg j ) } f_j p +(p)
\tilde{\partial}_j \\ 
& \equiv (p) \tilde{\partial}_j \mod ( \{ p \} \cup W, w).
\end{aligned}
\end{equation*} 
\end{Pf}

In the rest of this paper, we shall omit brackets whenever it is convenient.
Namely, the Lie product $[a, b]$ will be written as $ab$.  Moreover,
$(\mathrm{ad}x)^n y$ will be written as $x^n y$ and 
$x (\widetilde{\mathrm{ad}}y)^n$
as $xy^n$.  It would be clear from the context whether a product $ab$ 
means a Lie product or not. 

We write $f \equiv g \ \text {mod} \, (S, n)$ if $f-g=\sum \alpha_i
a_i s_i b_i$ with $l(a_i \overline{s}_i b_i)  \le n$, where
$n \in \mathbb{Z}_{>0}$, 
$\alpha_i \in k$, $a_i, b_i \in X^*$, and  $s_{i} \in S$.

\begin{Lem} \label{lem-b} Let $p \in S_+$.  
Then for any $l=1, \cdots, r$, 
we have $$(p)\tilde{\partial}_l \equiv 0 \quad \mod (S_+ \cup W, \,
l(\overline{p})).$$ \end{Lem}

\begin{Pf} 

\noindent
{\bf Case 1.} Relation $S_{+,1}$:

\vskip 2mm
Since $e_i^{1-n_{ij}}e_j
=\alpha e_je_i^{1-n_{ij}}$ with $\alpha \in k$, it suffices to prove our
assertion for $p = e_je_i^{1-n_{ij}}$ for $i \neq j$.  
We first consider the case when $a_{ii} = 2$.
We have only to check the cases when $l=i$ and $l=j$.  If $l=i$, we have 
\begin{equation*} 
\begin{aligned}
(p)\tilde{\partial}_i  = & (e_j e_i^{1-a_{ij}})
\tilde{\partial}_i \\ 
 =&  (e_j e_i^{ -a_{ij} }) h_i + (-1)^{\deg i} ( (e_j
e_i^{-a_{ij} -1}) h_i ) e_i \\ 
 &\mbox{ } + (-1)^{ 2 \deg i} ( (e_j
e_i^{-a_{ij} -2}) h_i ) e_i^2 + \dots + (-1)^{ -a_{ij} \deg i} (e_j h_i)
e_i^{-a_{ij}} \\ 
 \equiv & a_{ij} e_je_i^{-a_{ij} }+ (-1)^{\deg i}( a_{ij} -2)
e_je_i^{-a_{ij} } \\ 
& \mbox{ } + (-1)^{2 \deg i}( a_{ij} -4) e_je_i^{-a_{ij}
} + \dots + (-1)^{ -a_{ij} \deg i}( -a_{ij} ) e_je_i^{-a_{ij} }.
\end{aligned}
\end{equation*}
If $i \notin \tau$, 
then, clearly, the coefficient of $e_je_i^{-a_{ij} }$ is $0$.  
If $i \in \tau$, then $a_{ij} \in 2 \mathbb{Z}$ by the assumption on 
the generalized Cartan matrix $A$, 
and hence the coefficient of $e_je_i^{-a_{ij} }$ is also $0$.

Similarly, if $l=j$, we have 
\begin{equation*}
\begin{aligned}
(p)\tilde{\partial}_j & = (e_j e_i^{1-a_{ij}})
\tilde{\partial}_j  =  (-1)^{(1-a_{ij}) (\deg i) (\deg j) } h_j e_i^{1
-a_{ij} } \\ & \equiv  (-1)^{(1-a_{ij}) (\deg i) (\deg j) } a_{ji}e_i e_i^{
-a_{ij} } =0.
\end{aligned}
\end{equation*}

The proof for the case $a_{ii} =0$ is the same.

\vskip 2mm
\noindent
{\bf Case 2.} Relation $S_{+,2}$:

\vskip 2mm 
Let $p =(e_{k+1}e_k)(e_k e_{k-1})$
with $k \in \eta$.
If $l=k-1$, since $(e_{k+1}e_k)e_k $ or $e_{k+1}e_k$ is in $S_+$,
we have
\begin{equation*}
\begin{aligned}
(p)\tilde{\partial}_{k-1} & = ( (e_{k+1}e_k)(e_k e_{k-1})
)\tilde{\partial}_{k-1}  =  (e_{k+1}e_k)(e_k h_{k-1}) \\ 
& \equiv -a_{k-1, k}(e_{k+1}e_k)e_k \\ 
& \equiv  0 \ \  \mod \, (S_+ \cup W , \, l(\overline{p})).
\end{aligned}
\end{equation*}
Similarly, $(p) \tilde{\partial}_{k+1} \equiv 0 \ \ \mod \, 
(S_+ \cup W , \, l(\overline{p})).$

If $l=k$, since $a_{k,k-1} + a_{k,  k+1} = 0$ and $e_{k+1} e_{k-1} \in S_+$,
we have 
\begin{equation*}
\begin{aligned}
(p)\tilde{\partial}_{k} & = ( (e_{k+1}e_k)(e_k
e_{k-1}) )\tilde{\partial}_{k} \\ 
& =  (e_{k+1}e_k)(h_k e_{k-1}) -
(e_{k+1}h_k)(e_k e_{k-1}) \\ 
& \equiv a_{k,k-1}(e_{k+1}e_k)e_{k-1} + a_{k,k+1}
e_{k+1}(e_k e_{k-1}) \\ 
& =  (a_{k,k-1} + a_{k,  k+1}) e_{k+1}(e_k e_{k-1}) +
a_{k,k-1}(e_{k+1}e_k)e_{k-1} \\ 
& \equiv 0 \ \ \mod \, (S_+ \cup W , \, l(\overline{p})).
\end{aligned}
\end{equation*}
\end{Pf}

\begin{Lem} \label{lem-c} 
For any element $p \in S^c_+$ and $j =1, \cdots, r$, 
we have 
$$(p)\tilde{\partial}_j \equiv 0 \quad \mod (S^c_+ \cup W, \, 
l(\overline{p})).$$ 
\end{Lem}

\begin{Pf} 
As we have seen in Section 2, 
we have $S^c_+ = \bigcup S^{(i)}_+$ with
$S^{(i)}_+ \subset S^{(i+1)}_+$ for $i \ge 0$.  
Hence our assertion is equivalent
to saying that if $p \in S_{+}^{(i)}$, then 
$(p)\tilde{\partial}_j \equiv 0 \mod
(S^{(i)}_+ \cup W,  \, l(\overline{p}))$ for each $i \ge 0$.  
We will use induction on
$i$.  For $i=0$, it is simply Lemma \ref{lem-b}.  
Suppose that
$(q)\tilde{\partial}_j \equiv 0 \mod (S^{ (i) }_+ \cup W, 
l(\overline{q}) ) $ for
all $q \in S_{+}^{ (i) }$.  
Let $p \in S_{+}^{ (i+1) } \setminus S_{+}^{ (i) }$.  
Then
$ p = \langle q,r \rangle _{w} $ for some 
$q, r \in S_{+}^{ ( i ) }$
and $\langle q, r \rangle_{w} \equiv (q, r)_{w} 
\ \ \text {mod} \, (S^{ (i) }_+, w )$ 
by Lemma \ref{lem-equi}.
Since $l(w)=l(\overline{p})$,
we have 
$$\langle q, r \rangle_{w} \, \tilde \partial_{j} \equiv 
(q, r)_{w} \, \tilde\partial_{j} \ \ \text {mod} \
(S_{+}^{(i)} \cup W,  l(\overline {p})).$$
Thus  it is enough to
show that $(q, r)_{w} \, \tilde{\partial}_j \equiv 0 \mod(S^{(i)}_+
\cup W, l(\overline{p}))$.
Write $p=(q, r)_{w} = qa-br$.
Then by the induction hypothesis, we have 
\begin{equation*}
\begin{aligned} 
(q,r)_{w} \, \tilde{\partial}_j  = & q (a) \tilde{\partial}_j +
(-1)^{(\deg a) (\deg j)} (q)\tilde{\partial}_j a \\ 
 & - b (r)\tilde{\partial}_j - 
(-1)^{(\deg r) (\deg j)} (b)\tilde{\partial}_j r \\
& \equiv 0 \ \ \mod \, (S_{+}^{(i)} \cup W, \, l(\overline {p})).
\end{aligned}
\end{equation*}
\end{Pf}

Combining Lemma \ref{lem-a} and Lemma \ref{lem-c}, we obtain: 

\begin{Prop} \label{prop-a} 
For any element $p \in S^c_+$, we have 
$$\langle p, t_{ij} \rangle_{w} \equiv (p, t_{ij})_w \equiv 0 
\mod (S^c_+ \cup W, w).  $$
\end{Prop} 

Proposition \ref{prop-a} implies that all the compositions between the
relations in $S_{+}^c$ and $W$ are trivial. Similarly, one can show that 
all the compositions between the relations in $S_{-}^c$ and $W$ are 
also trivial. 
Now we can present the main theorem of this section.  

\begin{Thm} \label{thm-a} 
Let $\mathcal{G} = \mathcal{G}(A, \tau)$ be a Kac-Moody superalgebra
with the set of defining relations $S(A, \tau) = S_+ \cup W \cup S_-$.
Then the set $S^c_+ \cup W \cup S^c_-$ is a Gr\"obner-Shirshov
basis for the Kac-Moody superalgebra  $\mathcal{G}(A,\tau)$. 
That is, 
$S(A, \tau)^c = S^c_+ \cup W \cup S^c_-$. Hence it is also a Gr\"obner-
Shirshov basis for the universal enveloping algebra $\mathcal{U}(\mathcal{G})$
of $\mathcal{G}(A, \tau)$.
\end{Thm}

\begin{Pf} By definition, 
there is no nontrivial composition among the relations
in $S^c_{\pm}$ and the relations in $S^c_+$ and $S^c_-$.  
Also, all the compositions between the relations in $S_{\pm}^{c}$
and $W$ are trivial (see the remark after Proposition \ref{prop-a}).
Thus we have only to consider the compositions among the elements in $W$.  
We will show that $\langle p, q \rangle_{w} \equiv 0 \  
\text {mod} \, (W, w)$ for all $p, q \in W$, where $w\in X^*$ is determined 
by $p$ and $q$. 
There are four cases to be considered. 

If $p=h_i h_j \, (i>j)$ and $ q = h_j h_k \, (j>k)$, then 
$w= h_i h_j h_k$ and 
\begin{equation*}
\begin{aligned}
 \langle p, q \rangle _w & =  [w]_p -
[w]_q = (h_i h_j) h_k - h_i(h_j h_k) \\ 
& =  (h_i h_j) h_k \equiv 0.
\end{aligned}
\end{equation*}
If $p=e_j h_i + \alpha_{ij}e_j$ and $q=h_i h_k \, (i>k)$, then 
$w=e_j h_i h_k$ and 
\begin{equation*}
\begin{aligned}
\langle p, q \rangle _w & =  [w]_p - [w]_q =
(e_j h_i) h_k + a_{ij}e_j h_k - e_j(h_i h_k) \\ 
& =  (e_j h_k) h_i + a_{ij}e_j
h_k \equiv -a_{kj} e_j h_i + a_{ij}e_j h_k \\ 
& \equiv  a_{kj} a_{ij} e_j -
a_{kj} a_{ij}e_j =0.  
\end{aligned}
\end{equation*}
Similarly, if $p=h_i h_j \, (i>j)$ and $q= h_j f_k +a_{jk} f_k$, then 
$\langle p, q \rangle_{w} \equiv 0$. 
Finally, if $p=e_j h_i + a_{ij}e_j$ and  $q= h_i f_k + a_{ik} f_k$,
then $w=e_j h_i f_k$ and 
\begin{equation*} 
\begin{aligned}
\langle p, q \rangle _w & =  [w]_p - [w]_q 
= (e_j h_i) f_k + a_{ij}e_j f_k - e_j(h_i f_k) - a_{ik}e_j f_k\\ 
& =  (e_j f_k) h_i + a_{ij}e_j f_k -a_{ik} e_j f_k 
\\ & \equiv  \delta_{jk}h_j h_i + \delta_{jk} a_{ij} a_{ij} h_j 
- \delta_{jk}a_{ik} h_j \equiv 0, 
\end{aligned}
\end{equation*}
which completes the proof.  \end{Pf}

As a corollary, we obtain the {\it triangular decomposition} 
of Kac-Moody superalgebras and their universal enveloping algebras. 

\begin{Cor} \label{cor-a} 
Let $\mathcal{G}=\mathcal{G}(A, \tau)$ be a Kac-Moody superalgebra.
Then we have 
\begin{equation} \label {tri-a}
\mathcal{G} \cong \mathcal{G}_+ \oplus \mathcal{G}_0 \oplus
\mathcal{G}_-
\end{equation}
and 
\begin{equation} \label{tri-b}
U(\mathcal{G}) \cong U(\mathcal{G}_{+}) \otimes 
U(\mathcal{G}_{0}) \otimes U(\mathcal{G}_{-})
\end{equation}
as $k$-linear spaces. 
\end{Cor}

\begin{Pf} Observe that any super-Lyndon-Shirshov monomial of 
degree $\ge 2$ cannot be $W$-reduced if it contains
$h_i$ or $e_j f_k$ as a subword.  Hence by Theorem \ref{thm-a}, the
set $B$ of $S(A, \tau)^c$-reduced super-Lyndon-Shirshov monomials is 
given by $B=B_+ \cup H \cup B_-$,
where $B_{+}$ (resp. $B_{-}$) is the set of 
$S_{+}^c$-reduced (resp. $S_{-}^c$-reduced) super-Lyndon-Shirshov 
monomials in $e_i$'s (resp. $f_i$'s).  
By Theorem \ref{thm-b}, $B$ is a linear basis of
$\mathcal{G}$, which proves the $k$-linear isomorphism (\ref{tri-a}).
The isomorphism (\ref{tri-b}) follows from the Poincar\'e-Birkhoff-Witt
Theorem.  \end{Pf}

\section{Classical Lie superalgebras} \newcounter{r1} \setcounter{r1}{1}
\newcounter{r2} \setcounter{r2}{2} \newcounter{r3} \setcounter{r3}{3}
\newcounter{r4} \setcounter{r4}{4} \newcounter{r5} \setcounter{r5}{5}
\newcounter{r6} \setcounter{r6}{6} \newcounter{r7} \setcounter{r7}{7}
\newcounter{r8} \setcounter{r8}{8} \newcounter{r9} \setcounter{r9}{9}

In this section, we will give an explicit construction of 
Gr\"obner-Shirshov bases for the classical Lie superalgebras.
A Gr\"obner-Shirshov basis $S$ is said to be {\it minimal} if no proper
subset of $S$ is closed under the Lie composition. 
We first set up some notations.  Recall that we  omit
brackets whenever it is convenient. 
For the elements $x_{i} \in X$, we set
$[x_1 x_2 \dots x_m]= x_1[x_2 \dots x_m]$ and 
$\{x_1 \dots x_{m-1} x_m \} = \{x_1 \dots x_{m-1} \} x_m$ $(m\ge 1)$.
If $i>j$, we will write 
$x_{ij}=[x_{i} x_{i-1} \cdots x_{j}]$. 
For simplicity, we will also denote $x_{ii}=x_{i}$. 
We will use the lexicographical ordering for the set $I \times I$:
$(i,j)> (k,l)$ if and only if $i>k$ or $i=k$, $j>l$. 

We briefly recall the definition of classical Lie superalgebras 
\cite{Kac}.
Let $V=V_{\bar {0}} \oplus V_{\bar 1}$ be a $\mathbb {Z}_{2}$-graded 
vector space with $\dim V_{\bar 0}=m$ and $\dim V_{\bar 1}=n$, and let
$L$ be the space of $k$-linear endomorphisms of $V$.
For each $\alpha \in \mathbb{Z}_{2}$, set
$$L_{\alpha}=\{ T: V \rightarrow V | \ T(V_{\beta}) 
\subset V_{\alpha + \beta}  \ \ \text {for all} \ 
 \beta \in \mathbb{Z}_{2} \}.$$
Then $L$ has a $\mathbb{Z}_{2}$-graded decomposition
$L = L_{\bar{0}} \oplus L_{\bar{1}}$ and it becomes a Lie superalgebra
with the superbracket defined by 
$$[X, Y]=XY- (-1)^{\alpha \beta} YX$$ 
for $X \in L_{\alpha}$, $Y \in L_{\beta}$, $\alpha, \beta \in 
\mathbb{Z}_{2}$.
The Lie superalgebra $L$ is called the {\it general linear Lie 
superalgebra} and is denoted by $gl(m,n)$.

Let $v_{1}, \cdots, v_{m}$ be a basis of $V_{\bar 0}$ and $v_{m+1},
\cdots, v_{m+n}$ be a basis of $V_{\bar 1}$. Then $L$ can be interpreted 
as the space of $(m+n) \times (m+n)$ matrices over $k$, and we have 
$$ L_{\bar{0}} = \left\{ \left(\matrix A & 0 \\ 0 & D \endmatrix \right) 
\big| \ \text{$A$ is an $m \times m$ matrix and $D$ is 
an $n \times n$ matrix} \right\},$$
$$ L_{\bar{1}} = \left\{ \left( \matrix 0 & B \\ C & 0 \endmatrix \right) 
\big| \ \text {$B$ is an $m \times n$ matrix and $C$ is an 
$n \times m$ matrix} \right\}.$$

For $X = \left( \matrix A & B \\ C & D \endmatrix \right) \in gl(m,n)$, 
we define the {\it supertrace} of $X$ to be
$str X = tr A - tr B$, where $tr$ denotes the usual trace function.
Then the subspace $sl(m,n)$ of $gl(m,n)$ consisting of the matrices
with supertrace $0$ forms a Lie superalgebra which is
called the {\it special linear Lie superalgebra}. 

Let $B$ be a nondegenerate 
consistent supersymmetric bilinear form on $V$. Thus $V_{\bar 0}$
and $V_{\bar 1}$ are orthogonal to each other, 
$B|_{V_{\bar 0} \times V_{\bar 0}}$ is symmetric, and
$B|_{V_{\bar 1} \times V_{\bar 1}}$ is skew-symmetric (which implies
$n$ must be even). 
For each $\alpha \in \mathbb{Z}_{2}$, define
$$osp(m,n)_{\alpha}=\{ T \in gl(m,n)_{\alpha} | \
B(Tv, w)=-(-1)^{\alpha (\text {deg} v)} B(v, Tw) \ \ 
\text {for all} \ v, w \in V\}.$$
Then the subspace $osp(m,n)=osp(m,n)_{\bar 0} \oplus osp(m,n)_{\bar 1}$
becomes a Lie superalgebra. We set 
\begin{equation}
\begin{aligned}
B(m,n) &= osp(2m+1,2n) \ \ (m\ge 0, n>0), \\
C(n) & = osp(2, 2n-2) \ \ (n\ge 2), \\
D(m,n) &= osp(2m,2n) \ \ (m\ge 2, n>0).
\end{aligned}
\end{equation}
These subalgebras are called the 
{\it ortho-symplectic Lie superalgebras of type}
$B(m,n)$, $C(n)$, and $D(m,n)$, respectively.

\subsection{The special linear Lie superalgebra $sl(m,n) \, (m, n>0)$} 

Let $E_{ij}$ denote the $(m+n) \times (m+n)$ matrix whose
$(i, j)$-entry is equal to $1$ and all the other entries are $0$,
and let 
\begin{equation}
x_i = E_{i, i+1}, \ \ y_i = E_{i+1, i} \ \ \ 
(i=1,2,\cdots, m+n-1). 
\end{equation}
Then the elements $x_i $, $y_i$, $z_i=[x_i, y_i]$ $(i=1,2,\cdots, m+n-1)$ 
generate the Lie superalgebra $sl(m,n)$. 

On the other hand, let $\Omega = \{ 1, 2, \cdots , m+n-1 \}$,
$\tau = \{m \} \subset \Omega$, and consider the generalized Cartan matrix
$A=(a_{ij})_{i,j \in \Omega}$ defined by 
\begin{equation}
\begin{aligned}
\, & a_{m,m}=0, \ \ a_{m,m+1}=1, \ \ a_{m+1,m}=-1, \\
\, & a_{ij}=-1 \ \ \text {if} \ |i-j|=1 \ \text {and} \ (i,j) \neq (m,m+1), \\
\, & a_{ij}&=0 \ \ \text {if} \ |i-j|>1.
\end{aligned}
\end{equation}
Let $\mathcal{G}=\mathcal{G}(A, \tau)$ be the Kac-Moody superalgebra
associated with $(A, \tau)$ and denote by $e_i$, $f_i$, $h_i$ $(i=1, \cdots,
m+n-1)$ the generators of $\mathcal{G}$.
Then it is straightforward to verify that the generators
$x_i, y_i, z_i$  $(i=1, \cdots, m+n-1)$ of the Lie superalgebra
$sl(m,n)$ also satisfy the defining relations of the Kac-Moody algebra
$\mathcal{G} = \mathcal{G}(A, \tau)$. Hence there exists a surjective
Lie superalgebra homomorphism 
$\phi : \mathcal{G} \rightarrow sl(m,n)$ given by $e_i 
\mapsto x_i$, $f_i \mapsto y_i$, $h_i \mapsto z_i$ 
\,$(i=1,2,\cdots, m+n-1)$. 

In the following lemma, we will derive more 
``refined'' relations of $\mathcal{G}$,
which will be used to construct
a Gr\"obner-Shirshov basis for the special linear Lie superalgebra
$sl(m,n)$. 
Recall that we  use the notation 
$e_{ij}=[e_{i} e_{i-1} \cdots e_{j}]$ for $i>j$ and 
$e_{ii}=e_{i}$.

\begin{Lem} \label{lem-4aa} 
In the Kac-Moody superlagebra $\mathcal{G}=\mathcal{G}(A, \tau)$, we have
\begin{equation}
e_{ij} e_{kl} =\delta_{j-1, k} e_{il} \quad 
\text {for all} \  (i, j) \ge (k, l). \end{equation}
\end{Lem}

\begin{Pf} We will proceed in several steps.

\vskip 2mm
\noindent
{\bf Step 1:} For all $j>k+1$, we have $e_{ij} e_{kl}=0$. 

\vskip 2mm
By the Serre relations, we have $e_j e_l = 0$ for all $j>l+1$.
Next, fix $l$ and assume that $j>k+1$, $k>l$.
Then by the Jacobi identity and induction hypothesis, we get 
\begin{equation*}
e_{j}e_{kl}  =  e_{j} (e_k e_{k-1,l}) 
 =  (e_{j} e_{k})e_{k-1,l} + (-1)^d e_k (e_{j} e_{k-1, l})=0,
\end{equation*} 
where $d = (\text {deg} e_{j})(\text {deg} e_{k}) \in \mathbb{Z}_{2}$.
Finally, fix $j$ and assume that $i>j>k+1$.
Then the induction argument yields
\begin{equation*}
e_{ij}e_{kl}  =  (e_i e_{i-1,j}) e_{k l} 
 =  e_i (e_{i-1,j} e_{kl})+ (-1)^d (e_i e_{kl}) e_{i-1, j} 
=0,
\end{equation*} 
where $d = (\text {deg} e_{i})(\text {deg} e_{i-1}) \in \mathbb{Z}_{2}$.

\vskip 2mm
\noindent
{\bf Step 2:} For all $i,j,k \in \Omega$, we have 
$e_{ij} e_{j-1, k} =e_{ik}$.

\vskip 2mm
If $i=j$, there is nothing to prove.
If $i>j$, then by induction argument and Step 1,  we obtain
\begin{equation*}
\begin{aligned}
e_{ij} e_{j-1,k} & = 
(e_i e_{i-1, j}) e_{j-1,k}
=  e_i ( e_{i-1,j}\, e_{j-1,k}) + (-1)^d (e_i
e_{j-1, k}) e_{i-1,j} \\
&= e_{ii} e_{i-1,k} = e_{ik},
\end{aligned}
\end{equation*}
where $d = (\text {deg} e_{i})(\text {deg} e_{i-1}) \in \mathbb{Z}_{2}$.

\vskip 2mm
\noindent
{\bf Step 3:} For all $i>j$, we have 
$e_i e_{ij} =0$ and $e_{ij}e_j =0$. 

\vskip 2mm
By the Serre relations, we have 
$e_ie_{i, i-1} = 0$. 
If $i \ge j+2$, then Step 2 implies 
$e_{ij} = e_{i, i-1} e_{i-2, j}$.  
Hence by Step 1, we obtain 
\begin{equation*}
e_i e_{ij} =  e_i ( e_{i,i-1} e_{i-2,j}) 
 = (e_i e_{i,i-1}) e_{i-2,j} + (-1)^d e_{i,i-1}( e_i e_{i-2,j}) = 0, 
\end{equation*}
where $d = (\text {deg} e_{i})(\text {deg} e_{i,i-1}) \in \mathbb{Z}_{2}$.

Similarly, we get $e_{ij} e_{j}=0$ for $i>j$. 

\vskip 2mm
\noindent
{\bf Step 4:} For all $k, l\ge 1$, we have 
$h_i e_{i+k, i-l} =0$. 

\vskip 2mm
By the relations in $W$, we obtain 
$$h_i e_{i+k, i-l} = (a_{i, i+1} + a_{ii} + a_{i, i-1})
e_{i+k, i-l} = 0.$$

\vskip 2mm
\noindent
{\bf Step 5:} For all $i>j$, we have 
$e_{ij} e_{i-1} = 0$. 

\vskip 2mm
If $j=i-1$,
then by the Serre relations, we get $e_{ij} e_{i-1}=0$.
Suppose first that $j<i-1$ and
$i-1 \neq m$.  
The by Step 3, we obtain 
\begin{equation*}
\begin{aligned}
 (e_{ij} e_{i-1}) e_{i-1} & = ( (e_i e_{i-1, j}
)e_{i-1} ) e_{i-1} 
\\ & =  ( e_i ( e_{i-1, j} e_{i-1}) +(-1)^d (e_i e_{i-1})
e_{i-1, j} ) e_{i-1}\\ 
& =  (-1)^d (e_i e_{i-1} ) (e_{i-1, j} e_{i-1} ) +
(-1)^{d'} ((e_i e_{i-1} ) e_{i-1} ) e_{i-1, j} \\
&= 0,
\end{aligned}
\end{equation*}
where $d = (\text {deg} e_{i})(\text {deg} e_{i-1})$
and  $d' = (\text {deg} e_{i-1})(\text {deg} e_{i-1,j})$.
Multiplying both sides by $f_{i-1}$ yields
\begin{equation*} 
\begin{aligned}
0 & =  ( (e_{ij} e_{i-1} ) e_{i-1} ) f_{i-1}\\ 
& =  (e_{ij} e_{i-1} ) ( e_{i-1} f_{i-1} ) + ( ( e_{ij} e_{i-1} ) f_{i-1} )
e_{i-1}\\ 
& = ( e_{ij} e_{i-1} ) h_{i-1} + (e_{ij} h_{i-1} ) e_{i-1} +
( (e_{ij} f_{i-1} ) e_{i-1} ) e_{i-1}.
\end{aligned}
\end{equation*} 
The second summand
is equal to $0$ by Step 4.
Since $e_{ij} f_{i-1}$ is a scalar multiple of $e_i e_{i-2, j}$,
the third summand is also equal to $0$.  
By the Jacobi identity
and Step 4, the first summand yields $2 e_{ij} e_{i-1} = 0$, which 
proves our claim.  

If $j <i-1$ and $i-1 = m$, 
since $(e_{m+1}e_{m})(e_{m}e_{m-1})=0$ by the Serre relations, we get 
\begin{equation*}
\begin{aligned}
 e_{m+1, j} e_m &= e_{m+1} (e_{mj} e_m) - (e_{m+1} e_m) e_{mj}\\
&= -(e_{m+1}e_{m}) (e_{m}(e_{m-1} e_{m-2,j})) \\
 &=-(e_{m+1} e_m) ( (e_m e_{m-1} ) e_{m-2, j} ) 
\\ & = - ( (e_{m+1} e_m) (e_m e_{m-1} ) ) e_{m-2, j} + (
e_m e_{m-1} ) ( ( e_{m+1} e_m ) e_{m-2, j} )
& = 0. 
\end{aligned}
\end{equation*}

\vskip 2mm
\noindent
{\bf Step 6:} For all $n>k \ge 0$, $m>l \ge 0$, we have 
$e_{m+k, m-l} e_{m+k, m-l}=0$.

\vskip 2mm
Suppose  $k=0$.  If $l=0$,  then we have to show
that $e_m e_m =0$.  
Note that 
$$0= e_{m}(e_{m}e_{m-1}) = (e_m e_m) e_{m-1} - e_m
(e_m e_{m-1} ) = (e_m e_m ) e_{m-1}.$$ 
Multiplying both sides by $f_{m-1}$, we obtain
\begin{equation*}
\begin{aligned}
0 & = ( (e_m e_m ) e_{m-1} ) f_{m-1} \\ 
& =  (e_m e_m) (e_{m-1} f_{m-1} ) + ( (e_m e_m) f_{m-1} ) e_{m-1} \\ 
& =  (e_m e_m) h_{m-1} = 2e_m e_m,
\end{aligned}
\end{equation*}
which implies $e_m e_m =0$.  

Next, suppose $l>0$.  If
$e_{m, m-l} e_{m, m-l} = 0$, then
\begin{equation*}
\begin{aligned}
 0 &= ( ( e_{m, m-l}e_{m, m-l} ) e_{m-l-1} ) e_{m-l-1} \\ 
& = ( e_{m, m-l} e_{m, m-l-1} )e_{m-l-1}
+ ( e_{m, m-l-1} e_{m, m-l} ) e_{m-l-1} \\ 
& = 2 e_{m, m-l-1} e_{m, m-l-1},
\end{aligned}
\end{equation*}
which yields $e_{m, m-l-1} e_{m, m-l-1}=0$.
Hence, by the downward induction, we conclude 
$e_{m, m-l} e_{m, m-l} = 0$ for all $m>l \ge 0$.  

Finally, if $k>0$, then our assertion follows from 
the same downward induction argument as above.

\vskip 2mm
\noindent
{\bf Step 7:} For all $k \ge k'$, $l \le l'$, we have 
$e_{m+k, m-l} e_{m+k', m-l'} =0$.

\vskip 2mm
Suppose $k'=k$.  If $l=l'$, then our assertion was proved in Step 6. 
If $l<l'$ and $e_{m+k,m-l} e_{m+k, m-l'} =0$, then
\begin{equation*}
\begin{aligned}
0 & =  ( e_{m+k, m-l}e_{m+k, m-l'} ) e_{m-l' -1} \\ 
& =  e_{m+k, m-l} ( e_{m+k, m-l'-1}) + (e_{m+k,m-l} e_{m-l' -1} ) 
e_{m+k, m-l'} \\ 
& =  e_{m+k, m-l} e_{m+k, m-l' -1}.
\end{aligned}
\end{equation*}
Hence by the downward induction, we get 
$e_{m+k, m-l} e_{m+k, m-l'} =0$ for all $l \le l'$.  

If $k>k'$, our assertion follows by the same 
downward induction argument.

\vskip 2mm
\noindent
{\bf Step 8:} For all $i \ge j >1$, we have  $e_{ij} e_{i, j-1} =0$.

\vskip 2mm
  If $i=j$, then our assertion is just the Serre relation.
Suppose $i>j$ and $i+1 \neq m$.
Then if $e_{ij} e_{i, j-1} =0$, we have
\begin{equation*}
\begin{aligned}
 0 &= e_{i+1}( e_{i+1} ( e_{ij} e_{i, j-1} ) ) \\ 
&= e_{i+1}( e_{i+1, j} e_{i, j-1} ) + (-1)^d e_{i+1} ( e_{ij} e_{i+1, j-1} )\\
&= (-1)^{d'} e_{i+1, j} e_{i+1, j-1} + (-1)^d e_{i+1, j} e_{i+1, j-1} ,
\end{aligned}
\end{equation*}
where $d = (\deg e_{i+1} ) (\deg e_{ij})$ and 
$d' = (\deg e_{i+1} ) (\deg e_{i+1, j} )$.  Since $i+1 \neq m$,
 we have $e_{i+1, j} e_{i+1,j-1}=0$ and the induction argument 
gives our relations. 
If $i>j$, $i+1=m$ and $e_{ij} e_{i, j-1}=0$, then by Step 7, we get 
$e_{i+1,j} e_{i+1, j-1} = e_{mj} e_{m, j-1}=0$.
Hence our assertion follows from the induction. 

\vskip 2mm
\noindent
{\bf Step 9:} For all $k \neq j-1$, $(i, j) \ge (k, l)$, we have
 $e_{ij} e_{kl}=0$.

\vskip 2mm
Fix $k=i$.  If $l=j$, then our assertion holds by Step 6.  
If $l=j-1$,  then it is just Step 8.  
If $l <j-1$, then, by Step 1 and Step 8, we have 
\begin{equation*} 
\begin{aligned}
e_{ij} e_{il} & = e_{ij} (e_{i, j-1} e_{j-2, l} ) \\ 
& = ( e_{ij} e_{i, j-1} ) e_{j-2, l} + (-1)^d
e_{i,j-1} ( e_{ij} e_{j-2, l} ) =0,
\end{aligned}
\end{equation*}
where $d=(\deg e_{ij}) (\deg e_{i,j-1}) \in \mathbb{Z}_{2}$.

Suppose $k<i$.  If $j>k+1$, our assertion holds by Step 1. 
Let us assume $k \ge j$.  If $k=l$,  then we may assume
$k<i-1$ by Step 5, and we have
\begin{equation*}
\begin{aligned}
 e_{ij} e_k & =  ( e_{i, k+2} e_{k+1,j} ) e_k \\ 
& =  e_{i, k+2} ( e_{k+1, j} e_k ) + (-1)^d ( e_{i, k+2} e_k )
e_{k+1, j} =0.
\end{aligned}
\end{equation*}
We shall use induction on $k-l$.
Note that if $k>l$, then we have
$$ e_{ij} e_{kl}= e_{ij} ( e_k e_{k-1, l} ) = ( e_{ij}
e_k ) e_{k-1, l} + (-1)^d e_k ( e_{ij} e_{k-1,l} ),$$  
where $d=(\deg e_{ij}) (\deg e_{k}) \in \mathbb{Z}_{2}$.
The first summand is equal to $0$ by the case $k=l$.  
Consider the second summand.
If $j \neq k$, then it is 0 by the induction hypothesis. 
If $j =k$, then by Step 2, it is equal to
$$(-1)^d e_{k} (e_{ik} e_{k-1,l}) 
=(-1)^d e_{k} e_{il} =0.$$

\end{Pf}


\newcounter{ctr1} 
Let $X= E \cup H \cup F = \{ e_{i}, h_{i}, f_{i} | \ i\in \Omega \}$ be a 
$\mathbb{Z}_{2}$-graded set, where 
$\Omega=\{1,2, \cdots, m+n-1\}$ and $\tau=\{m\}$ is the set of odd
index.
Let $R_+$ be the set of relations in $E^{\#}$ given by:  

\ \ \ \ I.  \ \  $ e_i e_j \quad (i > j+1)$,

\ \ \ II. \ \ $ e_{ij} e_{i-1} \quad (i >j)$,

\ \ III. \ \ $ e_{ij} e_{i, j-1} \quad (i \ge j >1)$,

\ \ IV. \ \  $ e_{m+k, m-l} e_{m+k, m-l} \quad (n>k \ge 0, m > l \ge 0)$.

Let $R_{-}$ be the set of relations in $F^{\#}$ 
obtained by replacing $e_{ij}$'s in $R_{+}$ by $f_{ij}$', 
and let $R(A, \tau) = R_+\cup W \cup R_-$.
Consider the Lie superalgebra 
$L = \mathcal{L}_X / \langle R(A, \tau) \rangle$, where 
$\langle R(A, \tau) \rangle$ denotes the ideal in $\mathcal{L}_{X}$
generated by $R(A, \tau)$. 
Then, by Lemma \ref{lem-4aa}, there is a surjective Lie superalgebra 
homomorphism $\psi: L \rightarrow \mathcal{G}$ defined by
$e_{i} \mapsto e_{i}$, $h_{i} \mapsto h_{i}$, $f_{i} \mapsto f_{i}$
$(i\in \Omega)$. 
We now prove the main result of this subsection.

\begin{Thm}\label{thm-GB1}

The set $R(A, \tau)$ of relations in $\mathcal{L}_{X}$ is a  
Gr\"obner-Shirshov basis for the Lie superalgebra $L$.
\end{Thm}

\begin{Pf} Set $R=R(A, \tau)$.  As in the proof of Corollary \ref{cor-a}, the
set of $R(A, \tau)$-reduced super-Lyndon-Shirshov 
monomials is $B=B_+ \cup H \cup B_-$, where
$B_{\pm}$ is the set of $R_{\pm}$-reduced super-Lyndon-Shirshov 
monomials in $\mathcal{L}_E$
(resp.  in $\mathcal{L}_F$).  
We claim that the set of $R_+$-reduced 
Lyndon-Shirshov monomials in
$\mathcal{L}_E$ is 
$$B'_+ = \{ e_{ij} | m+n > i \ge j \ge 1 \}.$$ 
Let $w$ be an
$R_+$-reduced Lyndon-Shirshov monomial in $\mathcal{L}_{E}$.  
If $l(w)=1$, then there is nothing to prove.
Suppose that $l(w) >1$.  Then $w=uv$, where $u$, $v$ are 
$R_+$-reduced Lyndon-Shirshov monomials.
By induction, we have $w= e_{ij} e_{kl}$, where $i\ge j$,
$k\ge l$ and $(i, j) > (k, l)$ in the lexicographical ordering.  
Note that we must have $i>k$, for if $i=k$, 
then $j-1 \ge l$ and $\overline{e_{ij} e_{i,j-1}}$ is
a subword of $\overline{w}$.  
We will show that $k=j-1$ and $i=j$. 
If $k> j$, then $\overline{w}$ contains
$\overline{e_{k+1, j} e_k}$ as a subword, and if $k=j$, 
then $\overline{w}$ contains
$\overline{(e_{k+1} e_k) e_k}$ as a subword.  
Finally, if $k \le j-2$, then $\overline{w}$ contains
$\overline{e_j e_k}$ as a subword.  
Hence we must have $k=j-1$. 
Moreover, since $w$ is a 
Lyndon-Shirshov monomial, we must have  $i=j$. 
Therefore, we obtain $w=e_{il}$, which  proves our claim.

Now, let $w$ be an $R_+$-reduced super-Lyndon-Shirshov 
monomial in $\mathcal{L}_{E}$s.  Then $w$ is a Lyndon-Shirshov monomial or
$w=uu$ with $u$ a Lyndon-Shirshov monomial 
in $E^{\#}_{ \bar{1} }$.  If the latter is true, then, as we have 
seen in the previous paragraph, 
we have $u=e_{m+k, m-l}$ $(n>k \ge 0, m > l \ge 0)$,
in which case $w$ is not $R_+$-reduced by IV.  
Therefore we have $$B_+ = B'_+ = \{ e_{ij} |
m+n > i \ge j \ge 1 \}.$$ 
Similarly, we get $B_- = \{ f_{ij} | m+n > i \ge j \ge
1 \}$.  

By Lemma \ref{lem-aa}, $B$ spans $L$.  Since $\phi$ and $\psi$ are
surjective, we have $\mathrm{card}(B) \ge \dim sl(m,n)$.  But the number of
elements of $B$ is $(m+n)^2 -1$, which is equal to the dimension of $sl(m,n)$.
Thus $\phi$ and $\psi$ are isomorphisms and $B$ is a linear basis of $L$.  
Therefore, by
Proposition \ref{prop-tt}, $R$ is a Gr\"obner-Shirshov basis for $L$.  
\end{Pf}

\begin{Rmk} The proof of Theorem \ref{thm-GB1} shows that 
the Lie superalgebras $L$, $\mathcal{G}(A, \tau)$ and $sl(m,n)$ are
all isomorphic. Hence Theorem \ref{thm-GB1} gives a Gr\"obner-Shirshov
basis for the Lie superalgebra $sl(m,n)$. 
Our argument also shows that $R(A, \tau)$ is actually a minimal
Gr\"obner-Shirshov basis. 
\end{Rmk}

\subsection{The Lie superalgebras of
type $B(m,n)$ \, $(m , n>0)$}

Let $E_{ij}$ denotes the $(2m+2n+1) \times (2m+2n+1)$ matrix whose 
$(i,j)$-entry is 1 and all the other entries are 0. 
Set 
\begin{equation}
\begin{aligned}
\, & x_i  = E_{2m+i+1, 2m+i+2} -E_{2m+n+i+2, 2m+n+i+1} 
\ \  (1 \le i \le n-1), \\
\, & x_n  = E_{2m+n+1, 1} +E_{m+1, 2m+2n+1}, \\ 
\, & x_{n+i}  = E_{i, i+1} - E_{m+i+1, m+i} \ \  (1 \le i \le m-1), \\ 
\, & x_{m+n} = \sqrt{2} ( E_{m, 2m+1} - E_{2m+1, 2m} ) ,\\
\, & y_i = E_{2m+i+2, 2m+i+1} - E_{2m+n+i+1,2m+n+i+2} 
\ \ (1 \le i \le n-1), \\ 
\, & y_n = E_{1, 2m+n+1} - E_{2m+2n+1, m+1}, \\ 
\, & y_{n+i} = E_{i+1, i} - E_{m+i, m+i+1} \ \  (1 \le i \le m-1), \\
\, & y_{m+n} = \sqrt{2} ( E_{2m+1, m} - E_{2m, 2m+1} ).
\end{aligned}
\end{equation}
Then the elements $x_i$, $y_i$, $z_{i}=[x_{i}, y_{i}]$ $(i=1,2, \cdots, m+n)$
generate the ortho-symplectic  Lie superalgebra $B(m,n)=osp(2m+1,2n)$ 
$(m,n>0)$  and $x_{n}$, $y_{n}$ are the odd generators.

On the other hand, let 
$\Omega = \{1,2, \dots m+n\}$, $\tau = \{n \} \subset \Omega$,
and consider the generalized Cartan matrix 
$A=(a_{ i j } )_{i, j \in \Omega}$ defined by 
\begin{equation}
\begin{aligned}
\, & a_{n,n}  = 0, \quad a_{n, n+1} = 1, \quad a_{m+n, m+n-1} = -2, \\
\, & a_{ij}  =-1 \quad \text {if} \  | i-j |=1, \ (i,j)\neq (n,n+1),
\ (m+n, m+n-1), \\ 
\, & a_{ij}  =0 \quad \ \text {if} \  |i-j|>1. 
\end{aligned}
\end{equation}
Let $\mathcal{G}=\mathcal{G}(A, \tau)$ be the Kac-Moody superalgebra
associated with $(A, \tau)$ and denote by
$e_{i}$, $f_{i}$, $h_{i}$ $(i=1,2, \cdots, m+n)$ the generators 
of $\mathcal{G}$. Then, as in the case of $sl(m,n)$, one can verify
that the generators $x_{i}$, $y_{i}$, $z_{i}$ $(i=1,2,\cdots, m+n)$
of the Lie superalgebra $osp(2m+1,2n)$ satisfy the defining relations 
of the Kac-Moody superalgebra $\mathcal{G}(A,\tau)$.
Hence there exists a surjective  Lie superalgebra homomorphism
$\phi: \mathcal{G} \rightarrow osp(2m+1,2n)$ given by
$e_{i} \mapsto x_{i}$, $f_{i} \mapsto y_{i}$, $h_{i} \mapsto z_{i}$
$(i=1,2,\cdots, m+n)$. 
As in Section 4.1, we first derive more relations in $\mathcal{G}$,
which will be used to construct a Gr\"obner-Shirshov basis for the
ortho-symplectic Lie superalgebra $B(m,n)=osp(2m+1,2n)$ $(m,n>0)$.

 \setcounter{ctr1}{1} 
\begin{Lem} \label{lem-4cc} In the Kac-Moody superalgebra
$\mathcal{G}=\mathcal{G}(A,\tau)$, we have
\begin{equation} 
\begin{aligned}
\ & e_{ij} e_{kl}  = \delta_{j-1, k} e_{il} \ \ \text {if} \ 
(i, j) \ge (k, l), m+n>k, \\
\ & [ e_{m+n, i} e_{m+n, j} e_{m+n,k} ] =0 \ \ (i,j,k \in \Omega), \\
\ & ( e_{m+n, i} e_{m+n, j} ) ( e_{m+n, k} e_{m+n, l} ) =0
\ \ (i,j,k,l \in \Omega).
\end{aligned}
\end{equation}
\end{Lem}

\begin{Pf} 
As in Lemma \ref{lem-4aa}, we will prove our assertion in several steps. 
\vskip 2mm 
\noindent
{\bf Step 1:} For all $(i, j) > (k, l)$ and $m+n>k$, we have
 $e_{ij} e_{kl} = \delta_{j-1, k} e_{il}$. 

\vskip 2mm
If we remove the $(m+n)$-th row and the $(m+n)$-th column of $A$,
then we get the generalized Cartan matrix for the Lie superalgebra 
$sl(m,n)$. Thus we have only to consider the case when $i=m+n$.  
Suppose $k \le m+n-2$. If $j=m+n$, then $e_{m+n} e_{kl} =0$ 
as in Step 1 of the proof of Lemma \ref{lem-4aa}.
If $j < m+n$, then  we have
\begin{equation*}
\begin{aligned}
e_{m+n, j} e_{kl} &= ( e_{m+n} e_{m+n-1, j} ) e_{kl} \\
&= e_{m+n} ( e_{m+n-1, j} e_{kl} ) + (-1)^d ( e_{m+n} e_{kl} ) e_{m+n-1, j} \\
&= \delta_{j-1, k} e_{m+n} e_{m+n-1, l} = \delta_{j-1, k} e_{m+n, l},
\end{aligned}
\end{equation*}
where $d=(\deg e_{m+n}) (\deg e_{m+n-1,j})$.

If $k=m+n-1$ and $j=m+n$, then $e_{m+n} e_{m+n-1, l} = e_{m+n,l}$
and if $j< m+n-1$, then 
\begin{equation*}
\begin{aligned}
e_{m+n, j} e_{m+n-1, l} &=e_{m+n,j} ( e_{m+n-1} e_{m+n-2, l} )\\ 
&= ( e_{m+n, j} e_{m+n-1} ) e_{m+n-2, l}
+ (-1)^d e_{m+n-1} ( e_{m+n, j} e_{m+n-2, l} ) \\
&= ( e_{m+n, j} e_{m+n-1} ) e_{m+n-2, l} 
+ (-1)^d \delta_{j-1, m+n-2} e_{m+n-1} e_{m+n, l},
\end{aligned}
\end{equation*}
where $d=(\deg e_{m+n,j})(\deg e_{m+n-1})$. 
As in Step 5 of the proof of  Lemma \ref{lem-4aa}, 
we have $e_{m+n, j} e_{m+n-1} =0$, which proves our claim.

\vskip 2mm
\noindent
{\bf Step 2:} For all $i\in \Omega$, we have 
 $[e_{m+n} e_{m+n} e_{m+n, i} ] =0$.

\vskip 2mm
It is clear that $[e_{m+n} e_{m+n} e_{m+n}]=0$. 
Suppose that $[e_{m+n} e_{m+n} e_{m+n, i}]=0$ for $i <m+n$.  
Multiplying both sides by $e_{i-1}$, we obtain
\begin{equation*}
\begin{aligned} 
0 & =  [ e_{m+n} e_{m+n} e_{m+n, i} ] e_{i-1} \\
 &= e_{m+n} ( ( e_{m+n} e_{m+n, i} ) e_{i-1} ) + (-1)^d
(e_{m+n} e_{i-1} ) ( e_{m+n} e_{m+n, i} ) \\ 
&= e_{m+n} ( e_{m+n} e_{m+n, i-1}) 
+ (-1)^{d'} e_{m+n} ( ( e_{m+n} e_{i-1} ) e_{m+n, i} ) \\ 
& =  [ e_{m+n} e_{m+n} e_{m+n, i-1} ]
\end{aligned}
\end{equation*}
for $d, d' \in \mathbb{Z}_{2}$, and the downward induction on $i$  gives
our claim. 

\vskip 2mm
\noindent
{\bf Step 3:}  $[e_{m+n} e_{m+n} e_{m+n-1}](e_{m+n} e_{m+n-1} ) =0$.

\vskip 2mm If $m \neq 1$, then by the Serre relation,
we get 
\begin{equation*} 
\begin{aligned}
\, & [ e_{m+n} e_{m+n} e_{m+n-1}] ( e_{m+n} e_{m+n-1} ) \\ 
\, & = ( [e_{m+n} e_{m+n} e_{m+n-1} ]e_{m+n} ) e_{m+n-1} 
+ e_{m+n} ( [e_{m+n} e_{m+n} e_{m+n-1} ] e_{m+n-1} ) \\ 
\, & =  -[ e_{m+n} e_{m+n} e_{m+n} e_{m+n-1} ] e_{m+n-1} 
+ e_{m+n} ( e_{m+n} \{ e_{m+n} e_{m+n-1} e_{m+n-1} \} )=0. 
\end{aligned}
\end{equation*}
If $ m= 1$, then 
\begin{equation*}
\begin{aligned}  
\, & [e_{n+1} e_{n+1} e_n ] ( e_{n+1} e_n ) 
=  ( [e_{n+1} e_{n+1} e_n] e_{n+1} ) e_n 
+ e_{n+1} ( [ e_{n+1} e_{n+1}e_n ] e_n ) \\ 
\, & = e_{n+1} ( e_{n+1} \{ e_{n+1} e_n e_n \} ) 
- e_{n+1} ( (e_{n+1} e_n ) ( e_{n+1} e_n ) ) \\ 
\, & =  - [e_{n+1} e_{n+1} e_n ] ( e_{n+1} e_n )
- (e_{n+1} e_n ) [ e_{n+1} e_{n+1} e_n ] \\ 
\, &= -2 [e_{n+1} e_{n+1} e_n ] (e_{n+1} e_n ),
\end{aligned}
\end{equation*}
which yields $[e_{n+1} e_{n+1} e_{n}] [e_{n+1} e_{n}]=0$.

\vskip 2mm
\noindent
{\bf Step 4:} For all $i\in \Omega$, we have 
$$[e_{m+n} e_{m+n} e_{m+n-1} ] e_{m+n, i} =0, \ \ 
(e_{m+n} e_{m+n, i}) ( e_{m+n} e_{m+n-1})=0.$$

Let $a = [e_{m+n} e_{m+n} e_{m+n-1} ] e_{m+n, i}$ and $b = (e_{m+n}
e_{m+n, i}) ( e_{m+n} e_{m+n-1} )$.  
If $m \neq 1$, then by Step 2 and Step 3, we obtain
\begin{equation*} 
\begin{aligned}
0 & =  [e_{m+n} e_{m+n} e_{m+n, i} ] e_{m+n-1} \\ 
& =  e_{m+n} \{ e_{m+n} e_{m+n, i} e_{m+n-1} \} + ( e_{m+n}
e_{m+n-1} ) ( e_{m+n} e_{m+n, i} ) \\ 
& =  e_{m+n} \{ e_{m+n} e_{m+n-1}e_{m+m,i} \} - b = a - 2b,
\end{aligned}
\end{equation*}
and 
\begin{equation*} 
\begin{aligned}
0= &( [e_{m+n} e_{m+n}e_{m+n-1} ] ( e_{m+n} e_{m+n-1} ) ) e_{m+n-2, i} \\ 
 & =  [ e_{m+n} e_{m+n} e_{m+n-1} ] e_{m+n, i} 
+ ( [ e_{m+n} e_{m+n} e_{m+n-1} ] e_{m+n-2, i} ) (e_{m+n} e_{m+n-1} ) \\ 
 &=  a + e_{m+n} \{ e_{m+n} e_{m+n-1} e_{m+n-2, i} \} ( e_{m+n} e_{m+n-1} )
 = a + b .  
\end{aligned}
\end{equation*}
Hence we have $a=b=0$.  
Similarly, if $m=1$, then we get $a+b =0$ and $a+2b=0$,
which implies $a=b=0$.

\vskip 2mm
\noindent
{\bf Step 5:} For all $i, j \in \Omega$, we have 
$\{ e_{m+n} e_{m+n, i} e_{m+n, j} \} =0$.

\vskip 2mm
If $j=m+n$ or $j=m+n-1$, our assertion holds by Step 2 and Step 4.
We will use the downward induction on $j$. 
Suppose $j<m+n-1$ and $\{ e_{m+n} e_{m+n, i} e_{m+n, j} \} =0 $ for all
$i \in \Omega$.  
Then  we have
\begin{equation*}
\begin{aligned}
0 & =  \{ e_{m+n} e_{m+n, i} e_{m+n, j} \} e_{j-1} \\ 
& =  ( e_{m+n} e_{m+n, i} ) e_{m+n, j} + ( ( e_{m+n} e_{m+n, i} )
e_{j-1} ) e_{m+n, j} \\ 
& =  \{ e_{m+n} e_{m+n, i} e_{m+n, j-1} \} +
\delta_{ij} (e_{m+n} e_{m+n, i-1} ) e_{m+n, j} \\ 
& =  \{ e_{m+n} e_{m+n, i} e_{m+n, j-1} \}, 
\end{aligned}
\end{equation*}
which proves our claim. 

\vskip 2mm
\noindent
{\bf Step 6:} For all $i, j, k \in \Omega$, we have 
 $[ e_{m+n, i}e_{m+n, j} e_{m+n,k} ] =0$. 

\vskip 2mm
If $i=m+n$, Step 5 implies
$$ e_{m+n} ( e_{m+n, j} e_{m+n, k} ) =
(e _{m+n} e_{m+n, j} ) e_{m+n, k} + e_{m+n, j} ( e_{m+n} e_{m+n, k} ) =0.$$
If $i< m+n$, by the above observation, we get
\begin{equation*} 
\begin{aligned}
\, & [ e_{m+n, i} e_{m+n, j} e_{m+n, k} ]
= ( e_{m+n} e_{m+n-1, i} ) ( e_{m+n, j} e_{m+n, k} ) \\ 
\, & =  e_{m+n} [ e_{m+n-1,i} e_{m+n, j} e_{m+n, k} ] 
+ (-1)^d [ e_{m+n} e_{m+n, j} e_{m+n, k} ] e_{m+n-1,i} \\ 
\, & =  (-1)^{d'} \delta_{m+n, j} [ e_{m+n} e_{m+n, i} e_{m+n, k} ] +
(-1)^{d''} \delta_{m+n, k} [ e_{m+n} e_{m+n, j} e_{m+n, i}] = 0,
\end{aligned}
\end{equation*}
where $d, d', d'' \in \mathbb{Z}_{2}$.

\vskip 2mm
It remains to prove the last relation. But it is 
an immediate consequence of Step 6. 
\end{Pf}

Let $X=E\cup H\cup F =\{e_{i}, h_{i}, f_{i} | \ i\in \Omega \}$ be a
$\mathbb{Z}_{2}$-graded set, where $\Omega=\{1,2, \cdots, m+n \}$
and $\tau=\{n\}$ is the set of odd index. Let 
$R_+$ be the set of relations in $E^{\#}$ given by:  

\ \ \ \ I. \ \ $ e_i \, e_j \quad (m+n \ge i > j+1> 1)$, 

\ \ \ II. \ \ $ e_{ij} \, e_{i-1} \quad (m+n \ge i >j \ge 1) $,

\ \  III. \ \ $ e_{ij} \, e_{i,j-1} \quad (m+n > i \ge j >1) $,

\ \ IV. \ \ $ e_{n+k, n-l} \, e_{n+k, n-l} \quad (m>k \ge 0, \, n > l \ge 0)$,

\ \ \ V. \ \ $ [e_{m+n, i} \, e_{m+n, j} \, e_{m+n, j-1}] 
\quad (m+n \ge i \ge j >1)$,

\ \ VI. \ \  $ \{ e_{m+n, i} \, e_{m+n, j} \, e_{m+n, i-1}\} 
\quad (m+n \ge i >j \ge 1) $, 

\ VII. \ \  $ (e_{m+n, i} \, e_{m+n, j})( e_{m+n, i} \, e_{m+n, j}) 
\quad (m+n \ge i >n \ge j \ge 1) $, 

VIII. \ \ $ (e_{m+n, i} \, e_{m+n,j})( e_{m+n, i} \, e_{m+n, j-1 } )
 \quad (n \ge i > j > 1) $.

\vskip 2mm

Let $R_-$ be the set of relations in $F^{\#}$ obtained by replacing 
$e_{ij}$'s in $R_{+}$ by $f_{ij}$'s, and let $R(A, \tau) =
R_+ \cup W \cup R_-$.
Consider the Lie superalgebra $L = \mathcal{L}_X / \langle R(A, \tau)\rangle$.
Then, by Lemma \ref{lem-4cc}, there exists a surjective 
Lie superalgebra homomorphism $\psi: L \rightarrow \mathcal{G}$ 
defined by $e_{i} \mapsto e_{i}$, 
$h_{i} \mapsto h_{i}$ and $f_{i} \mapsto f_{i}$ $(i\in \Omega)$.
Then we have :

\begin{Thm} \label{thm-GB2} The set $R(A, \tau)$ of the relations in
$\mathcal{L}_{X}$ is a Gr\"obner-Shirshov basis for the
Lie superalgebra $L$.
\end{Thm}

\begin{Pf} Set $R= R(A, \tau)$.  
As in the case of $sl(m,n)$, the set of $R(A,
\tau)$-reduced super-Lyndon-Shirshov monomials 
is $B=B_+ \cup H \cup B_-$.  We claim that the set
of $R_+$-reduced Lyndon-Shirshov monomials 
in $\mathcal{L}_E$ is 
$$B'_+= \{ e_{ij} | i \ge j
\} \cup \{ e_{m+n, i} \, e_{m+n, j} | i>j \}.$$ 

Let $w$ be an $R_+$-reduced 
Lyndon-Shirshov
monomial in $\mathcal{L}_E$.  If $l(w) = 1$, there is nothing to prove.  
If $l(w)>1$, then $w=uv$, where $u$, $v$ are $R_+$-reduced 
Lyndon-Shirshov monomials. Hence  by
induction, we have  $u, v \in B'_+$.  
We will show that either $u=e_{m+n, j}$, $v=e_{m+n, l}$
with $j>l$ or $u=e_i$, $v= e_{i-1, l}$,
which would prove our claim.  
We need to consider the following four cases.

\vskip 2mm
\noindent
{\bf Case 1.}  $u = e_{ij}$, $v=e_{kl}$ \ \ $(i\ge j, k\ge l)$:

\vskip 2mm
Since $uv$ is Lyndon-Shirshov, we have
$(i,j) >(k, l)$ lexicographically.  If $i=k=m+n$,
then $u=e_{m+n, j}$, $v=e_{m+n, l}$ with $j>l$.  
If $i=k<m+n$, then $j-1 \ge l$, and $\overline{e_{ij} \,e_{kl}}$ 
contains $\overline{e_{ij} \, e_{i, j-1}}$ as a subword.
Hence $w$ is not $R_{+}$-reduced by III.  
If $i=j>k$ and $k=i-1$, then $u$ and $v$ have the desired form and 
we are done.
If $i=j>k$ and $k \le i-2$, then $w$ is not $R_{+}$-reduced by I.
If $i>k$ and $i>j$, then we must 
have $k=i-1$, since $e_{ij}=e_i \, e_{i-1, j}$ and $e_{i-1, j} \le e_{kl}$
by the definition of  Lyndon-Shirshov monomials.
Hence $w$ is not $R_{+}$-reduced by II.

\vskip 2mm
\noindent
{\bf Case 2.} $u=e_{kl}$, $v=e_{m+n, i} \, e_{m+n, j}$ 
\ \ $(i>j, k\ge l)$:

\vskip 2mm
Since $uv$ is Lyndon-Shirshov, we
have $k=m+n$ and $l \ge i$.  Then $w$ is not $R_{+}$-reduced by V, since
$\overline{w}$ contains $\overline{e_{m+n, i}
( e_{m+n, i} \, e_{m+n, i-1} ) }$ as a subword.

\vskip 2mm
\noindent
{\bf Case 3.} $u=e_{m+n, i} \, e_{m+n, j}$, $v=e_{kl}$ 
\ \ $(i>j, k\ge l)$: 

\vskip 2mm
Since $uv$ is Lyndon-Shirshov, we have
$e_{m+n, i} > e_{kl} \ge e_{m+n, j}$.  
It follows that $k=m+n$ and $i >l \ge j$.
Hence $\overline{w}$ contains 
$\overline{ \{e_{m+n, i} \, e_{m+n, j} \, e_{m+n, i-1}\} }$ as a 
subword, and $w$ is not $R_{+}$-reduced by VI.

\vskip 2mm
\noindent
{\bf Case 4.} $u=e_{m+n, i} \, e_{m+n, j}$, $v=e_{m+n, k} \, e_{m+n, l}$
\ \ $(i>j,  k>l)$: 

\vskip 2mm
Since $uv$ is Lyndon-Shirshov,  we have
$(i,j) > (k,l)$ and $e_{m+n, j} \le e_{m+n,k} \, e_{m+n, l}$. 
Thus  we have $j<k$ and either $i=k>j>l$ or $i>k>j$.  
If $i=k>j>l$,  then $\overline{w}$
contains $\overline{ (e_{m+n, i} \, e_{m+n, j} ) 
( e_{m+n, i} \, e_{m+n, j-1} )}$ as a subword, 
and $w$ is not $R_{+}$-reduced by VII or VIII.  
If $i>k>j$,  then
$\overline{w}$ contains 
$\overline{ \{ e_{m+n, i} \, e_{m+n, j} \, e_{m+n, i-1}\} }$ 
as a subword, and $w$ is not $R_{+}$-reduced by VI.

\vskip 2mm
Now, let $w$ be an $R_+$- reduced super-Lyndon-Shirshov monomial 
in $\mathcal{L}_{E}$.  Then $w$ is Lyndon-Shirshov
or $w=uu$ with $u$ a Lyndon-Shirshov monomial in $E^{\#}_{\bar{1} }$.  
If the latter is true, then we have the following three possibilities:

(i)  $u= e_{n+k, n-l}$ $(m>k \ge 0, n> l \ge 0)$,

(ii) $u=e_{m+n, j}$ $(1 \le j \le n)$, 

(iii) $u=e_{m+n, i} \, e_{m+n, j}$ $(m+n \ge i > n \ge j \ge 1)$.  

\noindent
But the cases (i) and (iii) cannot occur by IV and VII.
Therefore the set of $R_+$-reduced
super-Lyndon-Shirshov monomials is given by 
$$ B_+= \{e_{ij} | \ i\ge j \} \cup 
\{e_{m+n,i} e_{m+n,j}| \ i>j \} 
 \cup \{ e_{m+n, j} \, e_{m+n, j} | 1 \le
j \le n \}.$$ 
Similarly, we get 
$$ B_-= \{f_{ij} | \ i\ge j \} \cup 
\{f_{m+n,i} f_{m+n,j}| \ i>j \} 
 \cup \{ f_{m+n, j} \, f_{m+n, j} | 1 \le
j \le n \}.$$ 
By Lemma \ref{lem-aa} $B$ spans $L$.  
Since $\phi$ and $\psi$ are surjective, we have
$\mathrm{card}(B) \ge \dim osp(2m+1, 2n)$.  
But the number of elements of $B$ is
$2(m+n)^2 +m+3n$, which is equal to the dimension of $osp(2m+1,2n)$.  
Hence $B$ is a linear basis of $L$ and by Proposition \ref{prop-tt},
$R=R(A, \tau)$ is a Gr\"obner-Shirshov basis for the Lie superalgebra 
$L$. \end{Pf}

\begin{Rmk} The proof of Theorem \ref{thm-GB2} shows that 
the Lie superalgebras $L$, $\mathcal{G}(A, \tau)$ and 
$B(m,n)=osp(2m+1,2n)$ are
all isomorphic. Hence Theorem \ref{thm-GB2} gives a Gr\"obner-Shirshov
basis for the Lie superalgebra $B(m,n)=osp(2m+1,2n)$. 
Our argument also shows that $R(A, \tau)$ is actually a minimal
Gr\"obner-Shirshov basis. 
\end{Rmk}

\subsection{The Lie superalgebras of type $B(0,n) \, (n>0) $} 

Let $E_{ij}$ denotes the $(2n+1) \times (2n+1)$ matrix whose 
$(i,j)$-entry is 1 and all the other entries are 0.  
Set 
\begin{equation}
\begin{aligned}
x_i & = E_{i+1, i+2} - E_{n+i+2, n+i+1} \ \ (1 \le i \le n-1), \\ 
x_n &= \sqrt{2} (E_{1, 2n+1} + E_{n+1, 1} ),\\
y_i &= E_{i+2, i+1} - E_{n+i+1, n+i+2} \ \  (1 \le i \le n-1), \\ 
y_n &= \sqrt{2} (E_{1n} - E_{2n+1, 1} ).
\end{aligned}
\end{equation}
Then the elements $x_{i}$, $y_{i}$, $z_{i}=[x_{i}, y_{i}]$ 
generate the Lie superalgebra $B(0,n)=osp(1,2n)$ $(n>0)$ and $x_{n}$, $y_{n}$
are the odd generators.

On the other hand, let $\Omega= \{1,2, \dots n\} $, $\tau = \{n\} 
\subset \Omega$, and consider the generalized Cartan matrix 
$A=(a_{ i j } )_{i, j \in \Omega}$ defined by
\begin{equation}
\begin{aligned}
\, & a_{n,n} = 2, \quad a_{n, n-1} = -2,  \\ 
\, & a_{ij} =-1 \quad \text {if} \ | i-j |=1, \ (i,j) \neq (n, n-1), \\
\, & a_{ij}=0 \quad \ \text {if} \ |i-j|>1.  
\end{aligned}
\end{equation}
Let $\mathcal{G}=\mathcal{G}(A, \tau)$ be the
Kac-Moody superalgebra associated with $(A, \tau)$
and denote by $e_{i}$, $f_{i}$, $h_{i}$ $(i=1,2,\cdots, n)$ 
the generators of $\mathcal{G}$.
Then, by the same argument as in the proof of Lemma \ref{lem-4cc},
we obtain:

\begin{Lem} \label{lem-4dd} 
In the Kac-Moody superalgebra 
$\mathcal{G}=\mathcal{G}(A, \tau)$, we have
\begin{equation}
\begin{aligned}
\, & e_{ij}e_{kl} = \delta_{j-1, k} e_{il} 
\ \ \text {if} \ (i, j) \ge (k, l), n>k, \\
\, & [e_{ni} e_{n j} e_{n k} ] =0 \ \ (i,j,k\in \Omega), \\
\, & ( e_{n i} e_{n j} ) ( e_{n k} e_{n l} ) =0 
\ \ (i,j,k,l \in \Omega). 
\end{aligned}
\end{equation} \end{Lem}

Let $X=E\cup H \cup F =\{e_{i}, h_{i}, f_{i} | \ i\in \Omega \}$
be a $\mathbb{Z}_{2}$-graded set, where 
$\Omega=\{ 1, 2, \cdots, n \}$ and $\tau=\{n\} \subset \Omega$ is the
set of odd index. 
Let $R_+$ be the set of relations in $E^{\#}$  given by:  

\ \ \ \ I. \ \ $ e_i \, e_j \quad (n \ge i > j+1> 1)$, 

\ \ \ II. \ \ $ e_{ij} \, e_{i-1} \quad (n \ge i >j \ge 1) $,

\ \ III. \ \ $ e_{ij} \, e_{i,j-1} \quad (n > i \ge j >1) $,

\ \ IV. \ \ $ [e_{n, i} \, e_{n, j} \, e_{n, j-1}] \quad
(n \ge i \ge j >1)$,

\ \ \ V. \ \ $ \{ e_{n, i} \, e_{n, j} \, e_{n, i-1} \} 
\quad (n\ge i >j \ge 1) $, 

\ \ VI. \ \  $ (e_{n, i} \, e_{n, j})( e_{n, i} \, e_{n, j-1}) \quad
(n \ge i >j >1) $. 

Let $R_-$ be the set of relations in $F^{\#}$ obtained by 
replacing $e_{ij}$'s in $R_{+}$ by $f_{ij}$'s, 
and let $R(A, \tau) = R_+ \cup W \cup R_-$.
Consider the Lie superalgebra 
$L = \mathcal{L}_X / \langle R(A, \tau) \rangle$.  
Then there is a surjective Lie superalgebra homomorphism 
$\psi: L \rightarrow \mathcal{G}$, and using the same argument 
as in the proof of Theorem \ref{thm-GB2}, we obtain:

\begin{Thm} The set $R(A, \tau)$ of the relations in $\mathcal{L}_{X}$
is a Gr\"obner-Shirshov basis for the Lie superlagebra $L$.
\end{Thm}

\begin{Rmk}
The set of $R_{+}$-reduced super-Lyndon-Shirshov monomials in 
$\mathcal{L}_{E}$ is given by 
$$B_+= \{ e_{ij} | i \ge j \}
\cup \{ e_{n, i} \, e_{n, j} | i \ge j \},$$ 
and the Lie superalgebras $L$, $\mathcal{G}(A, \tau)$, and 
$B(0,n)=osp(1,2n)$ are all isomorphic. Moreover, $R(A, \tau)$
is a minimal Gr\"obner-Shirshov basis.
\end{Rmk}

\subsection{The Lie superalgebras of type $C(n) \, (n \ge 2)$} 

Let $E_{ij}$ denotes the $(2n+1) \times (2n+1)$ matrix whose 
$(i,j)$-entry is 1 and all the other entries are 0.  
Set 
\begin{equation}
\begin{aligned}
\, & x_1 = E_{13} - E_{n+2, 2}, \\
\, & x_i = E_{i+1, i+2} - E_{n+i+1, n+i} \quad (2 \le i \le n-1), \\
\, & x_n = E_{n+1, 2n}, \\
\, & y_1 = E_{31}+E_{2, n+2}, \\ 
\, & y_i = E_{i+2, i+1} - E_{n+i, n+i+1} \quad (2 \le i \le n-1), \\ 
\, & y_n = E_{2n, n+1}.
\end{aligned}
\end{equation}
Then the elements $x_{i}$, $y_{i}$, $z_{i}=[x_{i}, y_{i}]$
$(i=1,2,\cdots, n)$ are the generators of the Lie superalgebra
$C(n)=osp(2,2n-2)$, and $x_1, y_1$ are the odd generators.

Let $\Omega =\{1,2, \dots n\} $, $\tau = \{1 \} \subset \Omega$ and consider 
the generalized Cartan matrix $A=(a_{ i j } )_{i, j \in \Omega}$ 
defined by 
\begin{equation}
\begin{aligned}
\, & a_{11} = 0, \quad a_{1 2} = 1, \quad a_{n-1, n} = -2 , \\ 
\, & a_{ij} =-1 \quad \text {if} \ | i-j |=1, \ 
(i,j) \neq (1,2), (n-1,n), \\
\, & a_{ij} =0 \quad \ \text {if} \ |i-j|>1.
\end{aligned}
\end{equation}
Let $\mathcal{G}=\mathcal{G}(A, \tau)$ be the Kac-Moody superalgebra 
associated with $(A, \tau)$ and denote by $e_{i}$, $f_{i}$, $h_{i}$
$(i=1,2,\cdots, n)$ the generators of $\mathcal{G}$.
Then there is a surjective Lie superalgebra homomorphism 
$\phi: \mathcal{G} \rightarrow osp(2,2n-2)$ 
given by $e_{i} \mapsto x_{i}$, $f_{i} \mapsto y_{i}$,
$h_{i} \mapsto z_{i}$ $(i=1,2,\cdots, n)$. 

By a similar argument in the proof of Lemma \ref{lem-4cc},
we can derive a more refined set of relations in $\mathcal{G}$,
which gives a Gr\"obner-Shirshov basis for the Lie superalgebra 
$osp(2,2n-2)$ $(n\ge 2)$. 
Since the argument is a variation of the one given in 
Lemma \ref{lem-4cc}, we omit the proof here.

\begin{Lem} \label{lem-ee} 
In the Kac-Moody superalgebra $\mathcal{G}=\mathcal{G}(A, \tau)$, 
we have 
\begin{equation}
\begin{aligned}
\, & e_{ij} e_{kl} = \delta_{j-1, k} e_{il} 
\quad \text {if  $(i, j) \ge (k, l)$ and $k \neq n-1$ when $i=n$}, \\
\, & \{ e_{n i} e_{n-1, j} e_{n-1, k} \} =0 \ \ (n>i), \\
\, & \{ e_{ni} e_{n-1, j} e_{nk} \} =0 \ \ (i,j,k \in \Omega), \\
\, & ( e_{ni} e_{n-1, j} ) ( e_{n, k} e_{n-1, l} ) =0
\ \ (i,j,k,l \in \Omega).  
\end{aligned}
\end{equation}
\end{Lem}

Let $X=E\cup H \cup F =\{e_{i}, h_{i}, f_{i} | \ i\in \Omega \}$
be a $\mathbb{Z}_{2}$-graded set, where 
$\Omega=\{ 1, 2, \cdots, n \}$ and $\tau=\{1 \} \subset \Omega$ is the
set of odd index. 
Let $R_+$ be the set of relations in $E^{\#}$  given by:  

\ \ \ \ I. \ \ $ e_i \, e_j  \ \ (n \ge i > j+1> 1)$,

\ \ \ II. \ \ $ e_{ij} \, e_{i-1} \ \ (n > i >j \ge 1) $,

\ \ III. \ \ $ e_{ij} \, e_{i, j-1} \ \ (n \ge i \ge j >1) $,

\ \ IV. \ \ $ e_{i1} \, e_{i1} \ \ (n \ge i \ge 1) $,

\ \ \ V. \ \ $ \{ e_{n, i} \, e_{n-1, j} \, e_{n-1} \} \quad (n > j
\ge i \ge 1)$,

\ \ VI. \ \  $ \{ e_{n, i} \, e_{n-1, j} \, e_{n, i-1} \} \quad (n > j
\ge i > 1) $,

\ VII. \ \ $ [e_{n, i} \, e_{n, i} \, e_{n-1}] \quad (n \ge i> 1) $,

VIII. \ \ $ (e_{n, 1} \, e_{n-1, j})( e_{n, 1} \, e_{n-1, j } ) 
\quad (n > j > 1) $, 

\ \ IX. \ \ $ (e_{n, i} \, e_{n-1, j})( e_{n, i} \, e_{n-1, j-1 } ) 
\quad (n > j > i > 1) $.

\vskip 2mm
Let $R_-$ be the set of relations in $F^{\#}$ obtained by replacing 
$e_{ij}$'s in $R_{+}$ by $f_{ij}$'s, and 
let $R(A,\tau) = R_+ \cup W \cup R_-$.
Consider the Lie superalgebra 
$L = \mathcal{L}_X / \langle R(A, \tau) \rangle $.  
Then there is a surjective Lie superalgebra homomorphism 
$\psi: L \rightarrow \mathcal{G}$ defined by $e_{i} \mapsto e_{i}$,
$f_{i} \mapsto f_{i}$, $h_{i} \mapsto h_{i}$ $(i\in \Omega)$. 
Moreover, we have:

\begin{Thm} The set $R(A, \tau)$ of the relations in $\mathcal{L}_{X}$ 
is a Gr\"obner-Shirshov basis for the Lie superalgebra $L$.
\end{Thm}

\begin{Pf} Since our argument is similar to the one for the 
proof of Theorem \ref{thm-GB2}, we just give a sketch of the proof.
We first prove that the set of $R(A, \tau)$-reduced Lyndon-Shirshov 
monomials in $\mathcal{L}_{X}$ is given by
$$B'_+= \{ e_{ij} |\ i \ge j \} \cup \{ e_{n, i} \, e_{n-1, j}
| \ n >j \ge i \ge 1 \; \mbox{and}\; (i,j) \neq (1,1) \},$$
and conclude the set $B_{+}$ of $R_{+}$-reduced super-Lyndon-Shirshov 
monomials in $\mathcal{L}_{E}$ is equal to $B'_{+}$.

We see that $B=B_{+} \cup H \cup B_{-}$ spans $L$, where
$$B_-= \{ f_{ij} | \ i \ge j \} \cup \{ f_{n, i} \, f_{n-1, j}
| \ n >j \ge i \ge 1 \; \mbox{and}\; (i,j) \neq (1,1) \}$$
is the set of $R_{-}$-reduced super-Lyndon-Shirshov monomials in
$\mathcal{L}_{F}$. 
The number of elements in $B$ is $2n^2 +n-2$, 
which is equal to the dimension of
$osp(2,2n-2)$ $(n \ge 2)$.  
Hence the homomorphisms $\phi$ and $\psi$ are isomorphisms, and 
$B$ is a linear basis of $L$, which proves our assertion.  \end{Pf}

\begin{Rmk} The Lie superalgebras $L$, $\mathcal{G}(A, \tau)$, and 
$C(n)=osp(2,2n-2)$ are all isomorphic and $R(A, \tau)$
is a minimal Gr\"obner-Shirshov basis.
\end{Rmk}

\subsection{The Lie superalgebras of type $D(m,n) \, (m \ge 2, n>0) $}

Let $E_{ij}$ denotes the $(2m+2n) \times (2m+2n)$ matrix whose 
$(i,j)$-entry is 1 and all the other entries are 0.  
Set 
\begin{equation}
\begin{aligned}
\, & x_i = E_{2m+i, 2m+i+1} -E_{2m+n+i+1, 2m+n+i} \ \ 
 (1 \le i \le n-1), \\ 
\, & x_n = E_{2m+n, 1} + E_{m+1,2m+2n}, \\ 
\, & x_{n+i} = E_{i, i+1} - E_{m+i+1, m+i} \ \ (1 \le i \le m-1), \\
\, & x_{m+n} = E_{m, 2m-1} - E_{m-1, 2m}, \\
\, & y_i = E_{2m+i+1, 2m+i} - E_{2m+n+i, 2m+n+i+1} \ \  (1 \le i \le n-1), \\
\, & y_n = E_{1, 2m+n} - E_{2m+2n, m+1}, \\ 
\, & y_{n+i} = E_{i+1, i} -E_{m+i, m+i+1} \ \ (1 \le i \le m-1), \\ 
\, & y_{m+n} = E_{2m-1, m} - E_{2m, m-1}.
\end{aligned}
\end{equation}
Then the elements $x_{i}$, $y_{i}$, $z_{i}=[x_{i}, y_{i}]$
$(i=1,2,\cdots, m+n)$ are the generators of the Lie superalgebra
$D(m,n)=osp(2m,2n)$, and $x_n$, $y_n$ are the odd generators.

Let $\Omega = \{1,2, \dots m+n\} $, $\tau =\{n \}$, and 
consider the generalized Cartan matrix $A=(a_{ i j } )_{i, j \in \Omega}$
defined by
\begin{equation}
\begin{aligned}
\, & a_{nn} = 0, \quad a_{n,n+1} = 1, \quad
a_{m+n-2, m+n} = -1 , \\ 
\, & a_{m+n-1, m+n}= 0 \quad a_{m+n, m+n-2} = -1, \quad
a_{m+n, m+n-1} = 0, \\ 
\, & a_{ij} =-1 \quad \text {if} \  | i-j |=1, \  \ \text {and} 
\ (i,j) \neq (n,n+1),\\ 
\, & \qquad \qquad \qquad \qquad 
(m+n-1,m+n), (m+n, m+n-1), \\
\, & a_{ij} =0 \quad \ \text {if} \ |i-j|>1, \ \ \text {and} \ 
(i,j) \neq (m+n-2, m+n), \\
\, & \qquad \qquad \qquad \qquad (m+n, m+n-2).
\end{aligned}
\end{equation}

Let $\mathcal{G}=\mathcal{G}(A, \tau)$ be the Kac-Moody superalgebra 
associated with $(A, \tau)$ and denote by $e_{i}$, $f_{i}$, $h_{i}$
$(i=1,2,\cdots, m+n)$ the generators of $\mathcal{G}$.
Then there is a surjective Lie superalgebra homomorphism 
$\phi: \mathcal{G} \rightarrow osp(2m,2n)$ 
given by $e_{i} \mapsto x_{i}$, $f_{i} \mapsto y_{i}$,
$h_{i} \mapsto z_{i}$ $(i=1,2,\cdots, m+n)$. 

\vskip 2mm 
We modify some of our notations: 

(i) We neglect $e_{m+n, m+m-1}$; 
if $j \le m+n-2$, we write $e_{m+n, j}=
e_{m+n} \, e_{m+n-2, j}$.  

(ii) We introduce a modified Kronecker's delta:
$$\hat{\delta}_{ij} =\cases 
1 \ \  & \text {if} \ i=j \text {or} \ i=j+1=m+n-1,  \\ 
0 \ \  & \text {otherwise}.  \endcases.  $$

In the following lemma, we will list a set of relations in $\mathcal{G}$ 
which would yield a Gr\"obner-Shirshov basis for the Lie superalgebra
$D(m,n)=osp(2m,2n)$ $(m\ge 2, n>0)$. We will omit the proof which is 
similar to that of Lemma \ref{lem-4cc}.

\begin{Lem} \label{lem-4ff} 
In the Kac-Moody superalgebra $\mathcal{G}=\mathcal{G}(A, \tau)$, 
we have 
\begin{equation}
\begin{aligned}
\, & e_{ij} e_{kl} = \hat{\delta}_{j-1, k} e_{il} \ \ \text {if} 
\ (i, j) \ge (k, l), (i, k) \neq (m+n, m+n-1), \\
\, &  e_{m+n, i} e_{m+n-1, i} =0 \ \ \text {if} \ i > n, \\
\, & e_{m+n, i} e_{m+n-1, i-1} = e_{m+n, i-1} e_{m+n-1, i} \quad \text {if} 
\ i \le n, \\
\, & \{ e_{m+n, i} e_{m+n-1, j} e_{m+n-1, k} \} =0 \ \ (i,j,k\in \Omega), \\
\, & \{ e_{m+n, i} e_{m+n-1, j} e_{m+n, k} \} =0 \ \ (i,j,k \in \Omega), \\
\, & ( e_{m+n, i} e_{m+n-1, j} ) ( e_{m+n, k}
e_{m+n-1, l} ) =0 \ \ (i,j,k,l \in \Omega).  
\end{aligned}
\end{equation}
\end{Lem}

Let $X=E\cup H \cup F =\{e_{i}, h_{i}, f_{i} | \ i\in \Omega \}$
be a $\mathbb{Z}_{2}$-graded set, where 
$\Omega=\{ 1, 2, \cdots, n \}$ and $\tau=\{1 \} \subset \Omega$ is the
set of odd index. 
Let $R_+$ be the set of relations in $E^{\#}$  given by:  

\ \ \ \ I. \ \ $e_i \, e_j  \  (i > j+1, (i,j) \neq (m+n, m+n-2))$,\; 
$e_{m+n} \, e_{m+n-1}$, 

\ \ \ II. \ \ $ e_{ij} \, e_{i-1} \  (m+n > i >j ) $, \; $e_{m+n, j} \,
e_{m+n-2} \  (m+n-2 \ge j)$,

\ \ III. \ \  $ e_{ij} \, e_{i, j-1} \  (i \ge j >1) $
with $j \le m+n-2$ when $i=m+n$, 

\ \ \ \ \ \ \ \ \  $e_{m+n}(e_{m+n} \, e_{m+n-2})$, 

\ \ IV. \ \  $e_{m+n, i} \, e_{m+n-1, i} \  (m+n-2 \ge i >n) $,

\ \ \ \ \ \ \ \ \  
$e_{m+n, i} \, e_{m+n-1,i-1} - e_{m+n, i-1} \, e_{m+n-1, i} \ 
(i \le n)$, 

\ \ \ V. \ \ $ e_{n+k, n-l} \, e_{n+k, n-l} \  (m \ge k \ge 0, n> l \ge 0)$, 

\ \ VI. \ \ $ \{ e_{m+n, i} \, e_{m+n-1,
j} \, e_{m+n-1} \} \  (i < j < m+n) $, 

\ \ \ \ \ \ \ \ \   $\{ e_{m+n, i} \, e_{m+n-1, i} \,
e_{m+n-1} \} \  (i \le n )$, 

\ VII. \ \  $ \{ e_{m+n, i} \, e_{m+n, i} \, e_{m+n-1} \}
\  (m+n -2 \ge i)$, 

VIII. \ \  $ \{ e_{m+n, i} \, e_{m+n-1, j} \, e_{m+n, i-1} \} \ (1<i<j<m+n)$,

\ \ \ \ \ \ \ \ \  $ \{ e_{m+n, i} \, e_{m+n-1, i} \, e_{m+n, i-1} \}
 \  (i \le n)$, 

\ \ IX. \ \  $ (e_{m+n, i} \, e_{m+n-1, j} ) ( e_{m+n, i} 
\, e_{m+n-1, j-1 } ) \ (n+1<i+1 < j <m+n) $, 

\ \ \ \ \ \ \ \ \   $ (e_{m+n, i} \, e_{m+n-1, j})( e_{m+n, i} \, e_{m+n-1,
j } ) \  (i \le n , i<j )$.

\vskip 2mm

Let $R_-$ be the set of relations in $F^{\#}$ obtained by replacing 
$e_{ij}$'s by $f_{ij}$'s in $R_{+}$, and 
let $R(A,\tau) = R_+ \cup W \cup R_-$.
Consider the Lie superalgebra 
$L = \mathcal{L}_X / \langle R(A, \tau) \rangle $.  
Then there is a surjective Lie superalgebra homomorphism 
$\psi: L \rightarrow \mathcal{G}$ defined by $e_{i} \mapsto e_{i}$,
$f_{i} \mapsto f_{i}$, $h_{i} \mapsto h_{i}$ $(i\in \Omega)$, and we have:

\begin{Thm} The set $R(A, \tau)$ of the relations in $\mathcal{L}_{X}$ 
is a Gr\"obner-Shirshov basis for the Lie superalgebra $L$.
\end{Thm}

\begin{Pf} As in the case of $C(n)=osp(2,2n-2)$, we only give a brief 
sketch of the proof here. 
The set of $R_{\pm}$-reduced super-Lyndon-Shirshov monomials in
$\mathcal{L}_{E}$ (resp. $\mathcal{L}_{F}$) is given by
\begin{equation*}
\begin{aligned}
\, & B_+= \{ e_{ij} | \ i \ge j \} \cup \{ e_{m+n,i} \, e_{m+n-1,j} 
| \ \mbox{$i<j$ or $i=j \le n$ } \}, \\
\, & B_{-}=\{ f_{ij} | \ i \ge j \} \cup \{ f_{m+n,i} \, f_{m+n-1,j} 
| \ \mbox{$i<j$ or $i=j \le n$ } \}.
\end{aligned}
\end{equation*}
Hence the number of elements in the set of $R(A, \tau)$-reduced 
super-Lyndon-Shirshov monomials in $\mathcal{L}_{X}$ is 
$2(m+n)^2 -m+n$, which is equal to the dimension of the Lie superalgebra
$D(m,n)=osp(2m,2n)$. 
Therefore, $B$ is a linear basis of $L$ and $R(A, \tau)$ is 
a Gr\"obner-Shirshov basis for $L$.
\end{Pf}

\begin{Rmk} The Lie superalgebras $L$, $\mathcal{G}(A, \tau)$, and 
$D(m,n)=osp(2m,2n)$ are all isomorphic and $R(A, \tau)$
is a minimal Gr\"obner-Shirshov basis.
\end{Rmk}

\end{document}